\documentclass[11pt, reqno]{amsart}

\usepackage{amssymb,latexsym,amsmath,amsfonts}
\usepackage{amsmath,amscd,amsthm,amssymb}

\numberwithin{equation}{section}

\newcommand{\Zz}{{\mathbb Z}}
\newcommand{\Ff}{{\mathbb F}}
\newcommand{\Cc}{{\mathbb C}}

\newcommand{\p}{{\mathbf P}}
\newcommand{\q}{{\mathbf Q}}

\newcommand{\M}{{\mathbf M}}
\newcommand{\A}{{\mathbf A}}
\newcommand{\B}{{\mathbf B}}
\newcommand{\C}{{\mathbf C}}
\newcommand{\D}{{\mathbf D}}

\newcommand{\SO}{{\mathrm{SO}}}
\newcommand{\Orth}{\mathrm{O}}
 
\newcommand{\Om}{\mathrm{\Omega}}
 \newcommand{\SL}{{\mathrm{SL}}}
 \newcommand{\GL}{{\mathrm{GL}}}
 \newcommand{\PSL}{{\mathrm{PSL}}}
 \newcommand{\PGL}{{\mathrm{PGL}}}
\newcommand{\PG}{{\mathrm{PG}}}
\newcommand{\col}{{\mathrm{col}}}

\DeclareMathOperator{\Ker}{Ker}
\DeclareMathOperator{\Ima}{Im}
\DeclareMathOperator{\rank}{rank}

\newcommand{\Proof}{ \noindent{\bf Proof:}\quad }

\def\QED{\qed\medskip\par}
\newtheorem{Theorem} {Theorem} [section]
\newtheorem{Proposition} [Theorem] {Proposition}
\newtheorem{Lemma} [Theorem] {Lemma}
\newtheorem{Corollary} [Theorem] {Corollary}
\newtheorem{Conjecture}[Theorem]{Conjecture}

\newtheorem{Definition}[Theorem]{Definition}
\newtheorem{Remark} [Theorem] {Remark}

\def\Sq{\;\Box}
\def\Nsq{\;\not\!\!\Box}

\begin{document}

\title[Dimensions of Binary Codes]
{Proofs of Two Conjectures On the Dimensions of Binary Codes}

\author[Wu]{Junhua Wu }

\address{Department of Mathematics, Lane College, Jackson, TN 38301, USA}
\email{jwu@lanecollege.edu}


\begin{abstract}
Let $\mathcal{L}$ and $\mathcal{L}_0$ be the binary 
codes generated by the column $\Ff_2$-null spaces of 
the incidence matrices of external points versus 
passant lines and internal points versus secant lines 
with respect to a conic in 
$\PG(2,q)$, respectively. We confirm the conjectures 
on the dimensions of $\mathcal{L}$ and $\mathcal{L}_0$ 
using methods from both finite geometry and modular 
representation theory.
\end{abstract}

\maketitle

\section{Introduction}\label{intro}

Let $\Ff_q$ be the finite field of order $q$, where 
$q = p^e$, $p$ is a prime and $e\ge 1$ is an integer. 
Let $\PG(2,q)$ denote the classical projective plane 
of order $q$ represented via homogeneous coordinates. 
Namely, a point $\p$ of $\PG(2,q)$ can be written as
$(a_0,a_1,a_2)$, where $(a_0,a_1,a_2)$ is a non-zero 
vector of $V$, and a line $\ell$ as $[b_0, b_1, b_2]$,
where $b_0$, $b_1$, $b_2$ are not all zeros. The point 
$\p = (a_0, a_1, a_2)$ lies on the line 
$\ell=[b_0, b_1, b_2]$ if and only if
$$a_0b_0 + a_1b_1 + a_2b_2 = 0.$$

A non-degenerate conic in $\PG(2, q)$ is the set of 
points satisfying a non-degenerate quadratic form. 
It is well known that the set of points
\begin{equation}\label{conic}
\mathcal{O} = \{(1, t, t^2)\mid t \in \Ff_q\}\cup\{(0, 0, 1)\},
\end{equation}
which is also the set of projective solutions of the 
non-degenerate quadratic form
\begin{equation}\label{quad}
Q(X_0,X_1,X_2) = X_1^2 -X_0X_2 
\end{equation}
over $\Ff_q$, gives rise to a standard example of a non-degenerate 
conic in $\PG(2, q)$. It can be shown that every non-degenerate 
conic must has $q +1$ points and no three of them are collinear, 
which forms an oval (see \cite[P. 157]{hir}). In the case where 
$q$ is odd, Segre \cite{segre} proved that an oval in $\PG(2, q)$ 
must be a non-degenerate conic. In this paper, $q = p^e$ is 
always assumed to be an odd prime power. For convenience, we 
fix the conic in (\ref{conic}) as the ``standard" conic. A line 
$\ell$ is passant, tangent, or secant accordingly as 
$|\ell\cap \mathcal{O}|$ = $0$, $1$, or $2$, respectively. It is 
clear that every line of $\PG(2, q)$ must be in one of these 
classes. A point P is an internal, absolute, or external point 
depending on whether it lies on $0$, $1$, or $2$ tangent lines 
to $\mathcal{O}$. The sets of secant, tangent, and passant lines 
are denoted by $Se$, $T$ and $Pa$, respectively; the sets of 
external and internal points are denoted by $E$ and $I$, 
respectively. The sizes of these sets are 
$|Se| = |E| = \frac{q(q+1)}{2}$, $|Pa| = |I| = \frac{q(q-1)}{2}$, 
and $|T| = q+1$ (see (2.2)). Moreover, it can be shown that the 
quadratic form $Q$ in (\ref{conic}) induces a polarity $\sigma$, 
a correlation of order $2$, under which $E$ and $Se$, $O$ and $T$, 
and $I$ and $Pa$ are in one-to-one correspondence with each other, 
respectively. 

Let $\C$ be a 0-1 matrix; that is, $\C$ is a matrix whose entries 
are either $0$ or $1$. Note that $\C$ can be viewed as a matrix 
over any ring with $1$. The $p$-rank of $\C$, denoted by 
$\rank_p(\C)$, is defined to be the dimension of the column space 
of $\C$ over a field $F$ of characteristic $p$. The column null 
space of $\C$ over $F$ determines a linear code whose dimension 
is defined to be the dimension of the corresponding column null 
space of $\C$ over $F$.

Let $\A$ be the $(q^2 + q + 1) \times (q^2 + q + 1)$ point-line 
incidence matrix of $\PG(2, q)$; namely, $\A$ is a $0$-$1$ matrix 
and the rows and columns of $\A$ are labeled by the points and lines 
of $\PG(2, q)$, respectively, and the $(\p,\ell)$-entry of $\A$ is 
$1$ if and only if $\p\in\ell$. It can be shown that the $2$-rank 
of $\A$ is $q^2 + q$ \cite{HP} and the $p$-rank of $\A$ is
$\binom{p+1}{2}^e$ \cite{AK}, where $q = p^e$. The binary linear 
code generated by the column $\Ff_2$-null space of $\A$ has 
dimension $1$. Therefore, it is not useful for any practical 
purpose.

In \cite{keith1}, Droms, Mellinger and Meyer partitioned $\A$ 
into the following $9$ submatrices: 
\begin{equation}\label{partition}
\left(\begin{array}{ccc}
\A_{11}& \A_{12}& \A_{13}\\
\A_{21}& \A_{22}& \A_{23}\\
\A_{31}& \A_{32}& \A_{33}\\ 
\end{array}\right)
\end{equation}
where the block of rows for $(\A_{11},\A_{21},\A_{31})$ are 
labeled by the absolute, internal, and external points, 
respectively, and the block of columns for 
$(\A_{11}, \A_{12}, \A_{13})$ 
are labeled by the tangent, passant, and secant lines, respectively. 
They used the column null spaces of the submatrices $\A_{i,j}$ 
for $2\le i, j\le 3$ over $\Ff_2$ to construct four low-density 
parity-check (LDPC) codes. Based on computational evidence,
they made a conjecture on the dimensions of these codes.
For convenience, we denote $\A_{23}$ and $\A_{32}$ by $\B$ and 
$\B_0$, respectively. From (\ref{partition}), it follows that $\B$ 
and $\B_0$ are the incidence matrices of internal points versus
secant lines and external points versus passant lines, respectively. 
Note that $\B$ is a 
$\frac{q(q+1)}{2}\times\frac{q(q-1)}{2}$ matrix and $\B_0$ is 
a $\frac{q(q-1)}{2}\times\frac{q(q+1)}{2}$ matrix. The purpose of 
this article is to confirm the following conjecture on the 
dimensions of the LDPC codes $\mathcal{L}$ and 
$\mathcal{L}_0$ arising from the column $\Ff_2$-null spaces of 
$\B$ and $\B_0$, respectively.
\begin{Conjecture}
$($Droms, Mellinger and Meyer \cite{keith1}$)$\label{conj}
Let $\mathcal{L}$ and $\mathcal{L}_0$ be the $\Ff_2$-null 
spaces of $\B$ and $\B_0$, respectively. Then
\begin{displaymath}
\dim_{\Ff_2} (\mathcal{L}) =
\begin{cases}
\frac{q^2-1}{4}-q, & \text{if}\; q\equiv 1\pmod 4,\\
\frac{q^2-1}{4}-q+1, & \text{if}\; q\equiv 3\pmod 4;\\
\end{cases}
\end{displaymath}
and
\begin{displaymath}
\dim_{\Ff_2} (\mathcal{L}_0 ) =
\begin{cases}
\frac{q^2-1}{4}, & \text{if}\; q\equiv 1\pmod 4,\\
\frac{q^2-1}{4}+1, & \text{if}\; q\equiv 3\pmod 4.\\
\end{cases}
\end{displaymath}
\end{Conjecture}
Suppose that $\p_1$,..., $\p_{q(q+1)/2}$ 
and $\ell_1$,..., $\ell_{q(q-1)/2}$ are 
indexing the rows and columns of $\B$, 
respectively. Then we permute the rows 
and columns of $\B_0$ to obtain a new 
matrix $\C$ such that the rows and columns 
of $\C$ are indexed by $\ell_1$,..., 
$\ell_{q(q-1)/2}$ and $\p_1$,..., 
$\p_{q(q+1)/2}$, respectively. The matrix 
$\C$ is indeed equal to $\B^\top$, where 
$\B^\top$ is the transpose of $\B$. This 
implies that $\B$ and $\B_0$ have the same
$2$-rank. Therefore, in order to find the 
dimensions of the $\Ff_2$-null spaces of 
$\B$ and $\B_0$, it suffices to calculate 
the $2$-rank of either $\B$ or $\B_0$. 
Recall that the subgroup $G$ of $\PGL(3, q)$ 
fixing $\mathcal{O}$ is isomorphic to 
$\PGL(2, q)$ \cite[p. 158]{hir}. Further, 
$G$ has an index $2$ normal subgroup $H$ 
isomorphic to $\PSL(2, q)$. It is known 
\cite{dye} that $H$ acts transitively on 
$E$ and $I$ as well as on $Se$, $T$ and 
$Sk$.

In \cite{swx}, Sin, Wu and Xiang calculate 
the $2$-rank of $\A_{33}$ 
(i.e. the incidence matrix of external points 
and secant lines) using a combination of 
techniques from finite geometry and modular 
representation of $H$. In this article, 
we compute the $2$-rank of $\B$ using 
similar representation theoretic results 
obtained in \cite{swx} and different 
geometric results. Therefore, between 
the current article and \cite{swx}, 
the reader will expect to see some 
overlaps in the results and statements 
on modular representation of $H$ as 
well as the basic geometric facts 
about conics.

Let $F$ be an algebraic closure of $\Ff_2$. 
Let $F^I$ and $F^E$ be the free $F$-modules 
whose standard bases consist of the 
characteristic column vectors of $I$ and 
those of $E$, respectively. The actions of 
$H$ on $I$ and $E$ make the free $F$-modules 
$F^I$ and $F^E$ into $FH$-permutation 
modules. We define a map
\begin{equation}
\phi_{\B}: F^I\rightarrow F^E
\end{equation}
as follows: specify the images of the basis 
elements of $F^I$ under $\phi_{\B}$ first, i.e.
\begin{displaymath}
\phi_{\B}(\mathcal{G}_{\p}) =\displaystyle
\sum_{\q\in \p^\perp\cap E}\chi_{\q}
\end{displaymath}
for each $\p\in I$, and then extend $\phi_{\B}$ 
linearly to $F^I$, where $\perp$ is the polarity 
induced by the quadratic form $Q$, $\mathcal{G}_{\p}$ 
and $\chi_{\q}$ are the characteristic column 
vectors of the internal point $\p$ with respect
to $I$ and the external point $\q$ with respect
to $E$, respectively. The matrix of $\phi_\B$ 
is a $0$-$1$ matrix of size $|E|\times |I|$. 
Up to permutations of the rows and columns, 
$\B$ regarded as a matrix over $F$, is the 
matrix of $\phi_\B$ with respect to the 
standard bases of $F^I$ and $F^E$. Moreover, 
$\phi_{\B}(\mathbf{x}) = \B\mathbf{x}$ for 
$\mathbf{x}\in F^I$. It can be shown that 
$\phi_{\B}$ is an $FH$-homomorphism. Hence, 
the column space of $\B$ over $F$ is equal 
to $\Ima(\phi_{\B})$, which is also an 
$FH$-submodule of $F^E$. This point of view 
enables us to use results from modular 
representation of $H$ to determine the 
dimension of $\Ima(\phi_{\B})$ and thus 
the $2$-rank of $\B$. We remark that in 
the calculation of the $2$-rank of $\A_{33}$ 
the authors of \cite{swx} view $\A_{33}$ 
as the matrix of an $FH$-homomorphism $\phi$ 
from $F^E$ to $F^E$ under which the 
characteristic vector of an external point 
$\p$ is mapped to the sum of the characteristic 
vectors of the external points on $\p^\perp$.

Our idea of calculating 
$\dim_F (\Ima(\phi_{\B}))$ is to find a 
decomposition of $\Ima(\phi_{\B})$ into a 
direct sum of its submodules whose dimensions 
can be computed easily. To this end, we apply 
Brauer's theory and compute the decomposition
 of $\Ima(\phi_{\B})$ into blocks. 
The silmilar idea was used in \cite{swx} to 
compute the decomposition of $\Ker(\phi)$ 
into blocks as well as $\dim_F(\Ker(\phi))$. 
Nevertheless, there are two major differences 
between the current article and \cite{swx}: 
(1) the geometric results used to compute 
the decomposition of $\Ima(\phi_{\B})$ into 
blocks are essentially different from these 
used to compute the decomposition of 
$\Ker(\phi)$; (2) the summands of 
$\Ima(\phi_{\B})$ in its block decomposition 
are more complicated than these of $\Ker(\phi)$, 
which indicates that more efforts are required 
to find $\dim_F (\Ima(\phi_{\B}))$.

In the following we will give a brief 
overview of this article. In Section 2, 
we first review the basic facts about 
$\mathcal{O}$ and then prove several 
crucial geometric results. From them, 
in Section 5, we show that the $2$-rank 
of the incidence matrix $\D$ of external 
points and $N_{Pa,E}(\p)$ for $\p\in I$ 
(the set of external points on the 
passant lines through $\p$) is either 
$q$ or $q - 1$, depending on $q$. The 
character of the complex permutation 
module $\Cc^I$ and its decomposition 
into a sum of the irreducible ordinary 
characters of $H$ were calculated in 
\cite{wu}; the decomposition of the 
characters of $H$ into $2$-blocks was 
given by Burkhardt \cite{burkhardt} 
and Landrock \cite{landrock}. From 
them we see that $\Cc^I$ is a direct 
sum of $\Cc H$-modules consisting of 
one simple module from each block of 
defect zero, and some summands from 
blocks of positive defect. According 
to Brauer's theory, $\Ima(\phi_{\B})$ 
is the direct sum
\begin{equation}
\Ima(\phi_{\B}) =\displaystyle
\bigoplus_{B}\Ima(\phi_{\B})e_B
\end{equation}
where $e_B$ is a primitive idempotent 
in the center of $FH$. The block 
idempotents $e_B$ are elements of $FH$ 
and were computed in \cite{swx}. 
In order to compute $\Ima(\phi_{\B})e_B$ 
for each $2$-block $B$, we need detailed 
information concerning the action of group 
elements in various conjugacy classes on 
various geometric objects and on the 
intersections of certain special subsets 
of $H$ with various conjugacy classes 
of $H$. These computations are made in 
Sections 3 and 4. These information 
also tell us that 
(i) $\Ima(\phi_{\B})e_{B_0}$ is equal to 
the column space of $\D$ over $F$, or 
this space plus an additional trivial 
module, depending on $q$, where $B_0$ is 
the principal $2$-block of $H$, and (ii) 
block idempotents associated with 
non-principal $2$-blocks of positive 
defect annihilate $\Ima(\phi_{\B})$ 
(Lemma 6.2). Since the $B$-component 
of $F^I$ is the mod $2$ reduction of 
the $B$-component of $\Cc^I$, using 
(i) and (ii), and the block decomposition 
of $\Cc^I$ , we show that $\Ima(\phi_{\B})$ 
is equal to the direct sum of the column 
space of $\D$ and the simple modules 
lying in the $2$-blocks of defect $0$, 
or this sum plus an additional trivial 
module, depending on $q$. Then the 
dimension formula of $\Ima(\phi_{\B})$ 
follows instantly as a corollary.

\section{Geometric Results}

Recall that a {\it collineation} of 
$\PG(2,q)$ is an automorphism of 
$\PG(2,q)$, which is a bijection from 
the set of all points and all lines 
of $\PG(2,q)$ to itself that maps a 
point to a line and a line to a point, 
and preserves incidence. It is well 
known that each element of $\GL(3,q)$, 
the group of all $3\times 3$ 
non-singular matrices over $\Ff_q$, 
induces a collineation of $\PG(2,q)$. 
The proof of the following lemma is 
straightforward.

\begin{Lemma}\label{action1} Let
$\p=(a_0,a_1,a_2)$ 
$($respectively,$\;\ell$=$[b_0,b_1,b_2]$$)$ 
be a point 
$($respectively, a line$)$ of $\PG(2,q)$. 
Suppose that $\theta$ is a collineation 
of $\PG(2,q)$ that is induced by 
$\mathbf{D}\in GL(3,q)$. If we use 
$\p^\theta$ and $\ell^\theta$ to denote 
the images of $\p$ and $\ell$ under 
$\theta$, respectively, then
$$\p^\theta = (a_0,a_1,a_2)^\theta = (a_0,a_1,a_2)\mathbf{D}$$
and
$$\ell^\theta = [b_0,b_1,b_2]^\theta =[c_0,c_1,c_2],$$
where $c_0,c_1,c_2$ correspond to the 
first, the second, and the third coordinate 
of the vector 
$\mathbf{D}^{-1}(b_0,b_1,b_2)^\top$, 
respectively.
\end{Lemma}

A {\it correlation} of $\PG(2,q)$ is 
a bijection from the set of points to 
the set of lines as well as the set 
of lines to the set of points that 
reverses inclusion. A {\it polarity} 
of $\PG(2,q)$ is a correlation of 
order $2$. The image of a point $\p$ 
under a correlation $\sigma$ is 
denoted by $\p^\sigma$, and that of 
a line $\ell$ is denoted by $\ell^\sigma$. 
It can be shown \cite[p.~181]{hir} 
that the non-degenerate quadratic 
form $Q(X_0,X_1,X_2)$ = $X_1^2-X_0X_2$ 
induces a polarity 
$\sigma$ (or $\perp$) of $\PG(2,q)$, 
which can be represented by the matrix
\begin{equation}\label{matrix_M}
\M=\left(\begin{array}{ccc}
0 & 0 & -\frac{1}{2} \\
0 & 1 & 0 \\
-\frac{1}{2} & 0 & 0 \\
\end{array}\right).
\end{equation}

\begin{Lemma}\label{lemB}
$(${\rm{\cite[p.~47]{hp}}}$)$
Let $\p=(a_0,a_1,a_2)$ 
$($respectively, $\ell=[b_0,b_1,b_2])$ 
be a point $($respectively, a line$)$ 
of $\PG(2,q)$. If $\sigma$ is the polarity 
represented by the above non-singular 
symmetric matrix $\M$, then
$$\p^{\sigma} = (a_0,a_1,a_2)^{\sigma} =[c_0,c_1,c_2]$$
and
$$\ell^{\sigma} = [b_0,b_1,b_2]^{\sigma} =(b_0,b_1,b_2)\M^{-1},$$
where $c_0,c_1,c_2$ correspond to 
the first, the second, the third 
coordinate of the column vector 
$\M(a_0,a_1,a_2)^\top$, respectively.
\end{Lemma}
For example, if $\p=(x,y,z)$ is a point of $\PG(2,q)$,
then its image under $\sigma$ is $\p^{\sigma}=[z,-2y,x]$.

For convenience, we will denote the set 
of all non-zero squares of $\Ff_q$ by $\Sq_q$, 
and the set of non-squares by $\Nsq_q$. 
Also, $\Ff_q^*$ is the set of non-zero 
elements of $\Ff_q$.

\begin{Lemma}\label{bijection}
$(${\rm{\cite[p.~181--182]{hir}}}$)$
Assume that $q$ is odd.
\begin{itemize}
\item[(i)] The polarity $\sigma$ above 
defines three bijections; that is, 
$\sigma:\;I\rightarrow\;Pa$, $\sigma:\;E\rightarrow\;Se$, 
and 
$\sigma:\;\mathcal{O}\rightarrow\;T$ 
are all bijections.
\item[(ii)] A line $[b_0,b_1,b_2]$ of 
$\PG(2,q)$ is a passant, a tangent, or 
a secant to $\mathcal{O}$ if and only 
if $b_1^2-4b_0b_2 \in\Nsq_q$, 
$b_1^2-4b_0b_2 = 0$, or 
$b_1^2-4b_0b_2\in \Sq_q$, respectively.
\item[(iii)] A point $(a_0,a_1,a_2)$ 
of $\PG(2,q)$ is internal, absolute, 
or external if and only if 
$a_1^2-a_0a_2 \in \Nsq_q$, 
$a_1^2-a_0a_2 =0$, or 
$a_1^2-a_0a_2 \in \Sq_q$, respectively.
\end{itemize}
\end{Lemma}

The results in the following lemma 
can be obtained by simple counting; 
see {\rm \cite{hir}} for more details 
and related results.
\begin{Lemma}{\rm (\cite[p.~170]{hir})}
Using the above notation, we have
\begin{equation}\label{number}
|T|=|\mathcal{O}|=q+1,\;
|Pa|=|I|=\frac{q(q-1)}{2},\;\text{and}\;
|Se|=|E|=\frac{q(q+1)}{2}.
\end{equation}
Also, we have the following tables:
\begin{table}[htp]
\begin{center}
\caption{Number of points on lines of various types}
\bigskip
\begin{tabular}{cccc}
\hline
{Name}&{Absolute points}&{External points}&{Internal points} \\
\hline
{Tangent lines} & $1$ & $q$ & $0$ \\
{Secant lines}  & $2$ & $\frac{q-1}{2}$ & $\frac{q-1}{2}$ \\
{Passant lines} & $0$ & $\frac{q+1}{2}$ & $\frac{q+1}{2}$\\
\hline
\end{tabular}
\label{tab1}
\end{center}
\end{table}

\begin{table}[htp]
\begin{center}
\caption{Number of lines through points of various types}
\bigskip
\begin{tabular}{cccc}
\hline
{Name} & {Tangent lines} & {Secant lines} & {Skew lines} \\
\hline
{Absolute points} & $1$ & $q$ & $0$ \\
{External points} & $2$ & $\frac{q-1}{2}$ & $\frac{q-1}{2}$ \\
{Internal points} & $0$ & $\frac{q+1}{2}$ & $\frac{q+1}{2}$\\
\hline
\end{tabular}
\label{tab2}
\end{center}
\end{table}
\end{Lemma}

\subsection{More geometric results}

Let $G$ be the automorphism group of 
$\mathcal{O}$ in $\PGL(3,q)$
(i.e. the subgroup of $\PGL(3,q)$ 
fixing ${\mathcal O}$ setwise). 
Then $G$ is the image in $\PGL(3,q)$ 
of $\Orth(3,q)=\SO(3,q)\times\langle-1\rangle$,
hence also the image of $\SO(3,q)$, 
to which it is isomorphic. For our 
computations, we will describe $G$ 
in a slightly different way.
The map $\tau: \GL(2,q)\to \GL(3,q)$
sending the matrix 
$\left(\begin{smallmatrix}
a & b \\ 
c & d\end{smallmatrix}\right)$ to
\begin{equation}\label{3matrix}
\left(\begin{array}{ccc}
a^2 & ab & b^2\\
2ac & ad+bc & 2bd \\
c^2& cd & d^2\\
 \end{array}\right)
\end{equation}
is a group homomorphism. The image 
of $\tau(\GL(2,q))$ in $\PGL(3,q)$ 
lies in $G$. Now, whether or not 
the group $\tau(\GL(2,q))$ contains 
$\SO(3,q)$ depends on $q$. 
Nevertheless, $\tau(\GL(2,q))$ always 
contains a subgroup of index $2$ in 
$\Orth(3,q)$ whose image in $\PGL(3,q)$ 
is $G$. Thus, the induced homomorphism
$\overline\tau:\PGL(2,q)\to \PGL(3,q)$ 
maps $\PGL(2,q)$ isomorphically onto 
$G$.

Let $H=\tau(\SL(2,q))$, the group of 
matrices of the form (\ref{3matrix}) 
such that $ad-bc=1$. Since the kernel 
of $\tau$ is $\langle -I_2\rangle$, 
it follows that $H\cong \PSL(2,q)$ 
and that $H$ is isomorphic to its image 
$\overline{H}$ in $\PGL(3,q)$. In fact, 
we have $H=\Om(3,q)$. 
 
Since 
$$
\PGL(2,q)=\PSL(2,q)\cup 
\begin{pmatrix}1& 0\\ 0 & \xi^{-1}\end{pmatrix}\cdot \PSL(2,q),
$$
our discussion shows that
\begin{equation}\label{groupg}
H\cup {\bf d}(1,\xi^{-1},\xi^{-2})\cdot H
\end{equation}
is a full set of representative matrices 
for the elements of $G$. In our computations, 
it will often be convenient to refer to 
elements of $G$ by means of their 
representatives in the set (\ref{groupg}). 
Additionally, a group element in (\ref{3matrix}) 
has the inverse equal to
\begin{equation}\label{reverse}
\left(\begin{array}{ccc}
d^2 & -bd & b^2\\
-2cd & ad+bc & -2ab\\
c^2 & -ac & c^2\\
\end{array}\right).
\end{equation}
Moreover, the following holds.
\begin{Lemma}\label{transitive}\cite{dye}
The group $G$ acts transitively on $I$ and
$Pa$ as well as on $E$ and $Se$.
\end{Lemma}
We will refer to this lemma 
frequently in the rest of this section.

\begin{Lemma}\cite[Lemma 2.9]{swx}\label{meet}
Let $\p$ be a point not on ${\mathcal O}$, 
$\ell$ a non-tangent line, and $\p\in\ell$. 
Using the above notation, we have the
following.

\begin{enumerate}
\renewcommand{\labelenumi}{(\roman{enumi})}

\item If $\p\in I$ and $\ell\in$Pa, then 
$\p^\perp\cap \ell \in E$ if $q\equiv 1 \pmod 4$, 
and $\p^\perp \cap \ell \in I$ if $q\equiv 3
\pmod 4$.

\item If $\p\in I$ and $\ell \in$Se, then 
$\p^\perp\cap \ell \in I$ if $q\equiv 1 \pmod 4$, 
and $\p^\perp\cap \ell \in E$ if $q\equiv 3
\pmod 4$.

\item If $\p\in E$ and $\ell \in$Pa, then 
$\p^\perp\cap \ell \in I$ if $q\equiv 1 \pmod 4$, 
and $\p^\perp\cap \ell \in E$ if $q\equiv 3
\pmod 4$.

\item If $\p\in E$ and $\ell \in$Se, then 
$\p^\perp\cap \ell\in E$ if $q\equiv 1 \pmod 4$, 
and $\p^\perp\cap\ell\in I$ if $q \equiv 3
\pmod 4$.

\end{enumerate}
\end{Lemma}


Next we define $\Sq_q-1:=\{s-1\mid s\in \Sq_q\}$ and 
$\Nsq_q-1:=\{s-1\mid s\in \Nsq_q\}$.  

\begin{Lemma}\label{cs}\cite{store}
Using the above notation,
\begin{itemize}
\item[(i)] if $q\equiv 1\pmod 4$, then 
$|(\Sq_q-1)\cap\Sq_q|=\frac{q-5}{4}$ and $|(\Sq_q-1)\cap\Nsq_q|$ 
$=$ $|(\Nsq_q-1)\cap\Sq_q|$ $=$ $|(\Nsq_q-1)
\cap\Nsq_q|=\frac{q-1}{4}$;
\item[(ii)] if $q\equiv 3\pmod 4$, then $|(\Nsq_q-1)\cap\Sq_q|=\frac{q+1}{4}$ 
and $|(\Sq_q-1)\cap\Sq_q|$ $=$ $|(\Sq_q-1)\cap\Nsq_q|$ $=$ $|(\Nsq_q-1)\cap\Nsq_q|
=\frac{q-3}{4}$.
\end{itemize}
\end{Lemma}

\begin{Definition}
Let $\p$ be a point not on $\mathcal{O}$ and $\ell$ a line. 
We define $E_{\ell}$ $($respectively, $I_{\ell}$$)$ to be 
the set of external $($respectively, internal$)$ points on 
$\ell$, $Pa_{\p}$ $($respectively, $Se_{\p}$$)$ the 
set of passant $($respectively, secant$)$ lines through $\p$, 
and $T_{\p}$ the set of tangent lines through $\p$. Also, 
$N_{Pa, E}(\p)$ $($respectively, $N_{Se,E}(\p)$$)$ is defined 
to be the set of external points on the passant 
$($respectively, secant$)$ lines through $\p$.
\end{Definition}
In the following lemma, we list the sizes of the above
defined sets as well as the action of $G$ on these sets.
Also, we adopt standard notation from permutation group
theory. For instance, if $W\subseteq I$, then 
$W^g:=\{w^g\mid w\in W\}$, $G_\p$ is the stabilizer
of $\p$ in $G$, and for $M\subseteq G$, $M^g$ is the
conjugate of $M$ under $g$.
\begin{Lemma}\label{bsize}
Using the above notation, if $\p\in I$, we have
\begin{itemize}
\item[(i)] $|E_{\p^\perp}| = |Se_{\p}| = \frac{q+1}{2}$,
\item[(ii)] $|I_{\p^\perp}| = |Pa_{\p}| = \frac{q+1}{2}$, 
\item[(iii)] $|N_{Pa, E}(\p)| = |N_{Se, E}(\p)| = \frac{(q+1)^2}{4}$;
\end{itemize}
moreover, if $\p$ is not a point on $\mathcal{O}$, $\ell$
is a non-tangent line, and $g\in G$, we have
\begin{itemize}
\item[(iv)] $I_{\ell}^g=I_{\ell^g}$ and $Pa_{\p}^g=Pa_{\p^g}$, 
\item[(v)] $E_{\ell}^g=E_{\ell^g}$ and $Se_{\p}^g=Se_{\p^g}$, 
\item[(vi)] $H_\p^g=H_{\p^g}$, 
\item[(vii)] $N_{Pa, E}^g(\p) = N_{Pa, E}(\p^g)$ and 
$N_{Se, E}^g(\p)=N_{Se, E}(\p^g)$, 
\item[(viii)] $(\p^\perp)^g=(\p^g)^\perp$, where $\perp$ is the
polarity of $\PG(2,q)$ defined as above.
\end{itemize}

\end{Lemma}
{\Proof} The above (i) - (iii) follow from from Tables 1 and 2 and 
simple counting, and (iv) - (vii) follow from the fact that $G$
preserves incidence.
\QED

\noindent{By the defintion of $G$, it is clear that 
the following two lemmas are true.}
\begin{Lemma}
Let $\p$ be a point of $\PG(2, q)$. Then the polarity 
$\perp$ defines a bijection between $I_{\p^\perp}$ and 
$Pa_{\p}$, and also a bijection between $E_{\p^\perp}$ 
and $Se_{\p}$.
\end{Lemma}

\begin{Lemma}\label{stabilizer}
Let $W$ be a subgroup of $G$. Suppose that $g\in G$ 
and $\p$ is a point of $\PG(2,q)$. Then
$$(W^g)_{\p^g} = W_{\p}^g.$$
\end{Lemma}


\begin{Proposition}\label{Ktransitive}
Let $\p$ be a point not on $\mathcal{O}$ and set 
$K=G_\p$. Then $K$ is transitive on each of 
$I_{\p^\perp}$, $E_{\p^\perp}$, $Pa_{\p}$, and $Se_{\p}$. 
Moreover, if $\p\in E$, then $K$ is also transitive 
on $T_{\p}$.
\end{Proposition}
{\Proof} The case where $\p\in I$ is Proposition 2.11 
in \cite{wu}; the case where $\p\in E$ or $\mathcal{O}$ 
is Lemma 2.11(iii) in \cite{swx}.
\QED

\begin{Lemma}\label{a11}\cite[Corollary 2.16]{swx}
Let $\p$ be a point of $\PG(2,q)$ and let $\perp$ be 
the polarity of $\PG(2,q)$ defined above. Then for 
$g\in G_{\p}$ we have $\p^\perp=(\p^\perp)^g$. 
Consequently, $\p^\perp$ is fixed setwise by $G_{\p}$. 
Moreover, $G_{\p^\perp}=G_\p$.

\end{Lemma}

\begin{Lemma}\label{basic}
Assume that $\p\in I$ and $\ell=\p^\perp$. Let 
$\q\in E_{\ell}$ and $\ell^*\in T_{\q}$. Suppose 
that $\p_1$ and $\p_2$ are two distinct external 
points on $\ell^*$ and let $\ell_1$ and $\ell_2$ 
be the tangent lines different from $\ell^*$ 
through $\p_1$ and $\p_2$, respectively. Then 
$\ell_1$ and $\ell_2$ meet in an external point 
on a secant line through $\p$ if and only if one 
of the following two cases occurs:
\begin{itemize}
\item[(i)] $\p_1$ and $\p_2$ are on two passant 
lines through $\p$;
\item[(ii)] $\p_1$ and $\p_2$ are on two secant 
lines through $\p$.

\end{itemize}

\end{Lemma}

{\Proof} Since $G$ is transitive on $I$ and 
preserves incidence, without loss of generality, 
we may assume that $\p=(1,0,-\xi)$, and thus 
$\ell=[1,0,-\xi^{-1}]$. Since $K:=G_\p$ is 
transitive on $E_{\ell}$ by 
Proposition~\ref{Ktransitive}, we can assume that 
$\q=(0,1,0)$. Let $\ell^*=[1,0,0]$ be a tangent 
line through $\q$. It is clear that 
$$E_{\ell^*}=\{(0,1,m)\mid m\in\Ff_q\}.$$
Let $\p_1=(0,1,m_1)$ and $\p_2=(0,1,m_2)$ be two 
distinct external points on $\ell^*$. Then the 
tangent lines through $\p_1$ and $\p_2$ different 
from $\ell^*$ are $\ell_1=[m_1^2,-4m_1,4]$ and 
$\ell_2=[m_2^2,-4m_2, 4]$, respectively.  
So we have that  
$\p_3:=\ell_1\cap \ell_2=(1,\frac{m_1+m_2}{4}, 
\frac{m_1m_2}{4})\in E$.
Thus the line through $\p$ and $\p_3$ is
$$\ell_{\p,\p_3}=\left[m_1+m_2,-4\left(\frac{m_1m_2}
{4\xi}+1\right),\frac{m_1+m_2}{\xi}\right],$$
which is a secant line if and only if 
$$16\left(\frac{m_1m_2}{4\xi}+1\right)^2-
\frac{4(m_1+m_2)^2}{\xi}=\frac{(m_1^2-4\xi)
(m_2^2-4\xi)}{\xi^2}\in\Sq_q$$
if and only if either $m_i^2-4\xi\in\Nsq_q$ for $i=1$, 
$2$ or $m_i^2-4\xi\in\Sq_q$ for $i=1$, $2$. Since
the line through $\p$ and $\p_i$ ($i=1$ or $2$) is
$\ell_{\p,\p_i}=[1,-\frac{m_i}{\xi},\frac{1}{\xi}]$, 
and its discrimnant is $\frac{m_i^2-4\xi}{\xi^2}$,
we conclude that $\ell_{\p,\p_3}$ is a secant line 
if and only if either (i) $\p_1$ and $\p_2$ are on 
two passant lines through $\p$ or (ii) $\p_1$ and 
$\p_2$ are on two secant lines through $\p$.
\QED

\begin{Definition}
Let $N\subseteq E$. We define $\chi_N$ to be the 
characteristic $($column$)$ vector of $N$ with respect 
to $E$; that is, $\chi_N$ is a column vector of
length $|E|$ whose entries are indexed by the external 
points such that if $\p\in N$ then the entry of $\chi_N$ 
indexed by $\p$ is $1$, $0$ otherwise. For a line $\ell$, 
if no confusion occurs, we shoule use $\chi_{\ell}$ to 
replace $\ell_{E_{\ell}}$. Also, if $N=\{\p\}$ is a 
singleton set, we will frequently use $\chi_{\p}$ to 
replace $\chi_{\{\p\}}$.
\end{Definition}
\begin{Remark}
In the rest of this section, $\chi_N$ for $N\subseteq E$ 
will be always viewed as a column vector over $\Zz$, 
the ring of integer.
\end{Remark}

\begin{Corollary}\label{sksum}
Let $\p\in I$. Using the above notation, we have
$$\chi_{N_{Pa,E}(\p)}\equiv \displaystyle
\sum_{\ell\in T(\p,\ell(\p))}\chi_{\ell}\pmod 2,$$
where $\ell(\p)$ is a tangent line through an external 
point on $\p^\perp$, $T(\p,\ell(\p))$ is the set of 
tangent lines distinct from $\ell(\p)$ through the external 
points that are on both $\ell(\p)$ and the passant lines 
through $\p$, and the congruence means entrywise
congruence.

\end{Corollary}
{\Proof} It is clear that $|T(\p,\ell(\p))|=\frac{q+1}{2}$ 
since there are $\frac{q+1}{2}$ passant lines through $\p$ 
and each of them meets $\ell(\p)$ in an external point. Let 
$\ell\in T(\p,\ell(\p))$. Then by Lemma~\ref{basic}, any 
tangent line other than $\ell$ in $T(\p,\ell(\p))$ meets 
$\ell$ in an external point on a secant line through $\p$, 
and if we use $IE(\ell,\ell(\p))$ to denote their intersections 
with $\ell$ then the points in $E_{\ell}\setminus IE(\ell,\ell(\p))$
must be on the passant lines through $\p$. Since
$$(E_{\ell_1}\setminus IE(\ell_1,\ell(\p)))\cap (E_{\ell_2}
\setminus IE(\ell_2, \ell(\p)))=\emptyset$$ for two distinct 
lines $\ell_1$, $\ell_2\in T(\p,\ell(\p))$ and 
$$|E_{\ell}\setminus IE(\ell,\ell(\p))|=q-\frac{q-1}{2}=\frac{q+1}{2},$$
it follows that 
$$\displaystyle\sum_{\ell\in T(\p,\ell(\p))}|E_{\ell}\setminus 
IE(\ell,\ell(\p))|= 
\displaystyle\sum_{\ell\in T(\p,\ell(\p))}\frac{q+1}{2}
= \frac{(q+1)^2}{4}$$
which is the same as the size of $N_{Pa,E}(\p)$ by 
Lemma~\ref{bsize}(iii). Consequently, we must have
$$\displaystyle\bigcup_{\ell\in T(\p,\ell(\p))}E_{\ell}
\setminus IE(\p,\ell(\p)) = \displaystyle\bigcup_{\ell\in Pa_{\p}}
E_{\ell} = N_{Pa,E}(\p).$$
Moreover, since each point in $IE(\ell,\ell(\p))$ lies 
on exactly two lines in $T(\p,\ell(\p))$ and each point in 
$E_{\ell}\setminus IE(\ell,\ell(\p))$ doesn't lie on any 
line other than $\ell$ in $T(\p,\ell(\p))$, we obtain
\begin{equation}
\begin{array}{llllll}
\displaystyle\sum_{\ell\in T(\p,\ell(\p))} \chi_{E_{\ell}}& 
=  & \displaystyle\sum_{\ell\in T(\p,\ell(\p))}
\chi_{E_{\ell}\setminus IE(\ell,\ell(\p))}
+\displaystyle\sum_{\ell\in T(\p,\ell(\p))}
\displaystyle\sum_{\q\in IE(\ell,\ell(\p))} \chi_{\q}\\
{} & = & \displaystyle\sum_{\ell\in T(\p,\ell(\p))} 
\chi_{E_{\ell}\setminus IE(\ell, \ell(\p))} + 
2\displaystyle\sum_{\q\in M}\chi_\q\\
{} &\equiv & \displaystyle\sum_{\ell\in T(\p, \ell(\p))} 
\chi_{E_{\ell}\setminus IE(\ell, \ell(\p))}\\
{} & = & \displaystyle\sum_{\ell\in Pa_{\p}} \chi_{\ell} \\
{} & = & \chi_{N_{Pa, E}(\p)} \pmod 2

\end{array}
\end{equation}
where $M=\{\ell_1\cap \ell_2\mid \ell_1, 
\ell_2\in T(\p,\ell(\p)), \ell_1\not=\ell_2\}$.
\QED

\begin{Lemma}\label{set1}
Assume that $q\equiv 1\pmod 4$. Let $\p\in \mathcal{O}$. 
Then there exits a set $\mathcal{M}(\p)$ consisting of 
an odd number of internal points such that, for each 
external point $\q \in \p^\perp$, the number of passant 
lines through $\q$ and the points in $\mathcal{M}(\p)$, 
counted with multiplicity, is odd.
\end{Lemma}

\begin{Remark}
In this lemma, it is possible that $\q$, $\q_1$, ..., 
$\q_k$ are on the same passant line $\ell$, where 
$\q\in E_{\p^\perp}$ and $\q_i\in \mathcal{M}(\p)$ for 
$1\le i\le k$. If this circumstance occurs, then 
the line $\ell$ should be counted $k$ times.

\end{Remark}

{\Proof} Without loss of generality, we may assume 
that $\p=(0,0,1)$, and so $\ell:=\p^\perp=[1,0,0]$. 
Using Lemma~\ref{action1} and (\ref{reverse}), we 
have
\begin{equation}
K:=H_{\ell}=\left.\left\{\left(\begin{array}{ccc}
d^2& -bd& b^2\\
0 & 1 & -\frac{2b}{d} \\
0 & 0 & \frac{1}{d^2}\end{array}\right)
\right|d\in\Ff_q^*, b\in \Ff_q\right\}.
\end{equation}

Since
\begin{equation}
\left(\begin{array}{ccc} 
1 & -b & b^2\\ 
0 & 1 & -2b\\ 
0 & 0 & 1\end{array}\right)^k
=
\left(\begin{array}{ccc}
1 & -kb & (kb)^2\\ 
0 & 1 & -2kb \\ 
0 & 0 & 1 \end{array}\right)
\end{equation}
for any positive integer $k$, it is obvious that
\begin{equation}
\left.\left\{\left(\begin{array}{ccc}
1& -b & b^2\\ 
0 & 1 & -2b \\ 
0 & 0 & 1\end{array}\right)
\right| b\in \Ff_q\right\}
\end{equation}
is a collineation subgroup of order $q$ in $K$, 
which we denote by $T$.

For $(0,1,u_1)$, $(0,1,u_2)$ where $u_1$, $u_2\in\Ff_q$ 
and $u_1\not= u_2$, we have
\begin{equation*}
(0,1,u_1)\left(\begin{array}{ccc}
1 & -\frac{u_1-u_2}{2}& 
\left(\frac{u_1-u_2}{2}\right)^2\\
0 & 1 & -(u_1-u_2)\\
0 & 0 & 1 \end{array}\right)=(0,1,u_2);
\end{equation*}
this implies that $T$ is transitive on 
$E_{\ell}=\{(0,1,u)\mid u\in \Ff_q\}$.

Now let $\p_1=(1,0,-\xi)\in I$, set 
$\mathcal{M}(\p):=\{\p_1^g\mid g \in T\}$ 
which is the $T$-orbit of $\p_1$, and let 
$\q=(0,1,u)\in \ell$. Then
$$\mathcal{M}(\p)=\{(1,-b,b^2-\xi)\mid b\in \Ff_q\}$$
and the lines through both $\q$ and the points 
in $\mathcal{M}(\p)$ form the multiset
$$L(\q)=\{[b^2+ub-\xi,u,-1]\mid b\in \Ff_q\}.$$
Note that a line $[b^2+ub-\xi, u,-1]\in L(\q)$ is 
passant if and only if 
$\frac{(u+2b)^2}{4\xi}-1\in\Sq_q$.
Since the number of $t\in\Nsq_q$ satisfying 
$t-1\in \Sq_q$ is 
$|(\Nsq_q-1)\cap \Sq_q|=\frac{q-1}{4}$
by Lemma~\ref{cs}(i), it follows that the number 
of $b\in\Ff_q\setminus\{-\frac{u}{2}\}$ satisfying 
$\frac{(u+2b)^2}{4\xi}-1\in\Sq_q$ is 
$2(\frac{q-1}{4})=\frac{q-1}{2}$. Moreover, when 
$b=-\frac{u}{2}$, $\frac{(u+2b)^2}{4\xi}-1=-1\in\Sq_q$ 
as $q\equiv 1\pmod 4$. Hence, the number of $b\in\Ff_q$ 
satisfying $\frac{(u+2b)^2}{4\xi}-1\in\Sq_q$ is 
$\frac{q-1}{2}+1=\frac{q+1}{2}$. Thus, counted with 
multiplicity, there are $\frac{q+1}{2}$ passant lines 
in $L(\q)$. Therefore, there are an odd number of 
internal points (precisely $\frac{q+1}{2}$) in 
$\mathcal{M}(\p)$ connecting $\q$ by a passant
line as $q\equiv 1\pmod 4$. Since $T$ is transitive 
on both $\mathcal{M}(\p)$ and $E_{\ell}$ and 
preserves incidence, we conclude that the number 
of passant lines through an external point on 
$\p^\perp$ and the points in $\mathcal{M}(\p)$, 
counted with multiplictiy, must be odd.
\QED

\begin{Remark}
Let $\p \in \mathcal{O}$. In the rest of this 
article, without being further mentioned,
$\mathcal{M}(\p)$ always denotes a set of 
internal points associated with $\p$ satisfying 
the conditions in Lemma~\ref{set1}.
\end{Remark}

\begin{Corollary}\label{tsum1}
Assume that $q\equiv 1\pmod 4$. Let $\ell$ be 
a tangent line. Then
\begin{equation*}
\begin{array}{lll}
\chi_{\ell} & = & \displaystyle\sum_{\p\in\mathcal{M}(\ell^\perp)}
\displaystyle\sum_{\ell^{'}\in Pa_{\p}}\chi_{\ell^{'}}\\
{} & = & \displaystyle\sum_{\p\in\mathcal{M}(\ell^\perp)}
\chi_{N_{Pa, E}(\p)} 
\pmod 2
\end{array}
\end{equation*}
where the congruence is entrywise congruence.
\end{Corollary}

{\Proof} Let $\p\in \mathcal{M}(\ell^\perp)$. 
Then from Corollary~\ref{sksum},
it follows that
\begin{equation}
\begin{array}{lll}
\chi_{N_{Pa, E}(\p)} & = &\displaystyle
\sum_{\ell^{'}\in Pa_{\p}}\chi_{\ell^{'}}\\
{} & \equiv & \displaystyle\sum_{\ell^{'}
\in T(\p,\ell(\p))}\chi_{\ell^{'}} \pmod 2,
\end{array}
\end{equation}
where $\ell(\p)$ is a tangent line through 
an external point on $\p^\perp$ and 
$T(\p,\ell(\p))$ is the set of tangent lines 
different from $\ell(\p)$ through the 
external points that are both $\ell(\p)$ and 
the passant lines through $\p$.

Further, if we take $\ell(\p)=\ell$ for each 
$\p\in \mathcal{M}(\ell^\perp)$ and
set $W(\p):=\{\ell\cap\ell_1\mid \ell_1
\in Pa_{\p}\}$, then
\begin{equation}\label{eq1}
\begin{array}{lllllll}
\displaystyle\sum_{\p\in\mathcal{M}(\ell^\perp)}
\displaystyle\sum_{\ell^{'}\in Pa_{\p}}\chi_{\ell^{'}}
& \equiv & \displaystyle\sum_{\p\in\mathcal{M}(\ell^\perp)}
\displaystyle\sum_{\ell^{'}\in T(\p,\ell(\p))} \chi_{\ell^{'}}\\
{} & = & \displaystyle\sum_{\p\in\mathcal{M}(\ell^\perp)}
\displaystyle\sum_{\q\in W(\p)}
\displaystyle\sum_{\ell^{'}\in T_{\q}\setminus\{\ell\}}
\chi_{\ell^{'}}\\
{}& = & \displaystyle\sum_{\q\in E_{\ell}}
\displaystyle\sum_{\ell^{'}\in T_{\q}\setminus\{\ell\}}
a_{\ell^{'}}\chi_{\ell^{'}}\\
{} & = & \displaystyle\sum_{\ell^{'}\in T\setminus\{\ell\}}
a_{\ell^{'}}\chi_{\ell^{'}}\\
{} & \equiv & \displaystyle\sum_{\ell^{'}\in T\setminus\{\ell\}} 
\chi_{\ell^{'}} \pmod 2,
\end{array}
\end{equation}
where $a_{\ell^{'}}$ for $\ell^{'}\in T\setminus\{\ell\}$ 
are odd. In (\ref{eq1}), the second equality follows 
from the definition of $T(\ell, \ell(\p))$ and the 
third equality holds since the multiset
\begin{equation}\label{temp1}
\displaystyle\bigcup_{\q\in E_{\ell}}\displaystyle
\bigcup_{L(\q)}T_{\q}\setminus\{\ell\},
\end{equation}
where $L(\q):=\{\ell_{\p_1, \q}\in Pa \mid \p_1\in 
\mathcal{M}(\ell^\perp)\}$, is the same as the multiset
$$
\displaystyle\bigcup_{\p\in\mathcal{M}(\ell^\perp)}
\displaystyle\bigcup_{\q\in W(\p)}T_{\q}\setminus\{\ell\},
$$
and the tangent line $\ell^{'}$ other than $\ell$ 
through an external point $\q$ on $\ell$ occurs an 
odd number of times in (\ref{temp1}) by Lemma~\ref{set1}.

Since $\sum_{\ell^{'}\in T}\chi_{\ell^{'}}\equiv 0\pmod 2$,
it follows that
\begin{equation*}
\begin{array}{llllll}
\displaystyle\sum_{\p\in\mathcal{M}(\ell^\perp)}
\displaystyle\sum_{\ell^{'}\in Pa_{\p}}
\chi_{\ell^{'}}& \equiv & \displaystyle
\sum_{\p\in\mathcal{M}(\ell^\perp)}\chi_{N_{Pa, E}(\p)}\\
{} & \equiv & \displaystyle\sum_{\ell^{'}\in T
\setminus\{\ell\}}\chi_{\ell^{'}}\\
{} & \equiv & \chi_{\ell} + \displaystyle
\sum_{\ell^{'}\in T}\chi_{\ell^{'}}\\
{} & \equiv & \chi_{\ell} \pmod 2.
\end{array}
\end{equation*}
\QED

\begin{Lemma}\label{set22}
Let $\p\in E$ and let $T_1$ and $T_2$ be the two 
tangent lines through $\p$. Assume that 
$Z\subseteq (E_{T_1}\cup E_{T_2})\setminus \{\p\}$. 
Then there is a set $\mathcal{M}^{'}(\p)$ consisting 
of an even number of internal points such that, for 
any point $\q\in Z$, the number of passant lines 
through $\q$ and the points in $\mathcal{M}^{'}(\p)$, 
counted with multiplicity, is odd, and the number 
of passant lines through $\p$ and the points in 
$\mathcal{M}^{'}(\p)$, counted with mutiplicity, 
is even.
\end{Lemma}
{\Proof} Since $G$ is transitive on $E$, without 
loss of generality, we may assume that $\p=(0,1,0)$, 
and thus $T_1=[1,0,0]$ and $T_2=[0,0,1]$ are two 
tangent lines through $\p$. Let $K:=G_{\p}$ be 
the stabilizer of $\p$ in $G$. Using (\ref{groupg}), we have
\begin{equation}\label{stab}
\begin{array}{lllllll}
K & = & \left.\left\{{\bf d}\left(d^2, 1, \frac{1}{d^2}\right)
\right| d^2\in\Sq_q\right\}\cup
\left.\left\{{\bf ad}(\frac{1}{c^2}, -1, c^2)\right| 
c^2\in \Sq_q\right\}\\
{}&\cup&\left.\left\{{\bf d}\left(d^2, \frac{1}{\xi}, \frac{1}{d^2\xi^2}\right)\right| d^2\in\Sq_q\right\}\cup
\left.\left\{{\bf ad}(\frac{1}{c^2}, -\frac{1}{\xi}, \frac{c^2}{\xi^2})\right| c^2\in \Sq_q\right\}.
\end{array}
\end{equation}

Let $\p_1=(1,1,x)$, where $x\in \Nsq_q$ 
(respectively, $x\in \Sq_q$) and $1-x\in\Nsq_q$, 
be an internal point for $q\equiv 3\pmod 4$ 
(respectively, $q\equiv 1\pmod 4$). 
(Note that such an $x$ in the last coordinate 
of $\p_1$ exists in $\Ff_q$.) Then the $K$-orbit 
of $\p_1$ is
$$\mathcal{O}_{\p_1}=\left.\left\{
\left(1,\frac{1}{d^2}, \frac{x}{d^4}\right)\right| 
d^2\in\Sq_q\right\}
\cup\left.\left\{\left(1,\frac{1}{\xi d^2}, 
\frac{x}{\xi^2 d^4}\right)\right| d^2\in \Sq_q\right\}.$$

To prove the first part of  the lemma, we need 
only show that it holds for 
$$Z=(E_{T_1}\cup E_{T_2})\setminus\{\p\}.$$ Let 
$\q=(0,1,1)\in Z$. Using (\ref{stab}), we have that 
$K_{\q}$ only contains the identity collineation. 
So $K$ is transitive on $Z$ as $|Z|=|K|=2(q-1)$. 
The lines through $\q$ and the points in 
$\mathcal{O}_{\p_1}$ form the multiset
$$L(\q)=\{[x-d^2, d^4, -d^4]\mid d^2\in\Sq_q\}
\cup\{[x-d^2\xi, d^4\xi^2, -d^4\xi^2]\mid d^2\in\Sq_q\}.$$
A line in $L(\q)$ is passant if and only if 
$$\frac{(d^2-2)^2}{4(1-x)}-1\in\Sq_q$$
or
$$\frac{(d^2\xi-2)^2}{4(1-x)}-1\in\Sq_q,$$
where $d^2\in\Sq_q$. The number of $d^2$ satisfying 
either of the above two equations is equal to that 
of $t\in\Ff_q^*$ satisfying $\frac{(t-2)^2}{4(1-x)}-1\in\Sq_q$ 
since $\Ff_q^*=\Sq_q\cup \Sq_q\xi$, where 
$\Sq_q\xi=\{d^2\xi\mid d^2\in \Sq_q\}$. 
For the case where $q\equiv 3\pmod 4$, since the number 
of $t\in \Ff_q$ satisfying $\frac{(t-2)^2}{4(1-x)}-1\in\Sq_q$ 
is equal to 
$2|(\Nsq_q-1)\cap \Sq_q|=2(\frac{q+1}{4})=\frac{q+1}{2}$
by Lemma~\ref{bsize}(ii) and $t=0$ is one of them, 
we see that the number of passant lines in $L(\q)$, 
counted with multiplicity, is $\frac{q-1}{2}$ which is odd 
since $q\equiv 3\pmod 4$. For the case where 
$q\equiv 1\pmod 4$, since the number of 
$t\in\Ff_q\setminus\{2\}$ satisfying 
$\frac{(t-2)^2}{4(1-x)}-1\in \Sq_q$ is equal to 
$2|(\Nsq_q-1)\cap \Sq_q|=2(\frac{q-1}{4})=\frac{q-1}{2}$ 
by Lemma~\ref{bsize}(i), 
$t=0$ is not one of the solutions and $t=2$ 
also satisfies $\frac{(t-2)^2}{4(1-x)}-1\in\Sq_q$, 
we see that the number of passant lines in $L(\q)$, 
counted with multiplicity, is $\frac{q+1}{2}$, 
which is odd as $q\equiv 1\pmod 4$. Now we set 
$\mathcal{M}^{'}(\p):=\mathcal{O}_{\p_1}$,
and so $|\mathcal{M}^{'}(\p)|=q-1$ is even. 
Since $K$ is transitive on both 
$Z=(E_{T_1}\cup E_{T_2})\setminus\{\p\}$ and 
the points in $\mathcal{M}^{'}(\p)$, the number
of passant lines through a point in $Z$ and 
the points in $M^{'}(\p)$, counted with 
multiplicity, must be odd.

The lines through $\p$ and the points in 
$\mathcal{M}^{'}(\p)$ form the multiset
$$\left.\left\{\left[1,0,-\frac{d^4}{x}\right]\right| d^2\in\Sq_q\right\}\cup
\left.\left\{\left[1,0,-\frac{d^4\xi^2}{x}\right]\right| d^2\in\Sq_q\right\},$$
each or none of which is a passant line 
accordingly as $q\equiv 3\pmod 4$ or
$q\equiv 1\pmod 4$. Hence, we conclude 
that the number of passant lines through 
$\p$ and the points in $\mathcal{M}(\p)$, 
counted with multiplicty, is even.
\QED

\begin{Remark}
Let $\p\in E$. In the following discussion, 
without being further mentioned, 
$\mathcal{M}^{'}(\p)$ will always denote a 
set consisting of an even number of internal 
points that satisfy the conditions with 
$Z=E_{T_1}\setminus\{\p\}$ in the above lemma, 
where $T_1$ is one of the two tangent lines 
through $\p$.
\end{Remark}

\begin{Corollary}\label{tsum2}
Let $\p \in E$ and let $T_1$ and $T_2$ be 
the two tangent lines through $\p$. Then
\begin{equation*}
\begin{array}{lll}
\chi_{T_1}+\chi_{T_2} &\equiv & \displaystyle\sum_{\q\in\mathcal{M}^{'}(\p)}
\displaystyle\sum_{\ell\in Se_{\q}}\chi_{\ell}\\
{} &\equiv & \displaystyle\sum_{\q\in\mathcal{M}^{'}(\p)}\chi_{N_{Se, E}(\q)}\pmod 2,
\end{array}
\end{equation*}
where the congruence means entrywise congruence.
\end{Corollary}
{\Proof} Let $\q\in\mathcal{M}^{'}(\p)$. Then 
Corollary~\ref{sksum} gives
\begin{equation}
\begin{array}{lll}
\chi_{N_{Pa, E}(\q)}& \equiv & \displaystyle\sum_{\ell^{'}\in Pa_{\q}}\chi_{\ell^{'}}\\
{} & \equiv & \displaystyle\sum_{\ell^{'}\in T(\q, \ell(\q))} \chi_{\ell^{'}} \pmod 2,
\end{array}
\end{equation}
where $\ell(\q)$ is a tangent line through an 
external point on $\q^\perp$ and $T(\q,\ell(\q))$ 
is the set tangent lines through the external 
points that are on both $\ell(\q)$ and the passant 
lines through $\q$. Let ${\bf 1}$ be the all-one 
column vector of length $|E|$. Since
\begin{equation}
{\bf 1} +\chi_{N_{Pa, E}(\q)} \equiv \chi_{N_{Se, E}(\q)} \pmod 2
\end{equation}
and $|\mathcal{M}^{'}(\p)|$ is even, we have
\begin{equation}\label{temp3}
\begin{array}{llllll}
\displaystyle\sum_{\q\in\mathcal{M}^{'}(\p)}\displaystyle\sum_{\ell\in Se_{\q}}\chi_{\ell}
& \equiv & \displaystyle\sum_{\q \in \mathcal{M}^{'}(\p)}({\bf 1}+\chi_{N_{Pa, E}(\q)})\\
{} & \equiv & \displaystyle\sum_{\q \in \mathcal{M}^{'}(\p)}{\bf 1}+
\displaystyle\sum_{\q\in \mathcal{M}^{'}(\p)}\chi_{N_{Pa, E}(\q)}\\
{} &\equiv & \displaystyle\sum_{\q\in \mathcal{M}^{'}(\p)}\chi_{N_{Pa,E}(\q)}\\
{}&\equiv &\displaystyle\sum_{\q\in \mathcal{M}^{'}(\p)}
\displaystyle\sum_{\ell\in T(\q,\ell(\q))}\chi_{\ell}\pmod 2.
\end{array}
\end{equation}
Further, if we set $\ell(\q):=T_1$ for each 
$\q\in\mathcal{M}^{'}(\p)$ and set
$W^{'}(\q):=\{T_1\cap\ell_1\mid \ell_1\in Pa_{\q}\}$, 
since the multiset
\begin{equation*}
\displaystyle\bigcup_{\p_1\in E_{T_1}}\displaystyle\bigcup_{L^{'}(\p_1)}
T_{\p_1}\setminus\{T_1\},
\end{equation*}
where $L^{'}(\p_1)=\{\ell_{\p_1,\p_2}\in Pa\mid \p_2\in\mathcal{M}^{'}(\p)\}$,
is the same as the multiset
\begin{equation}\label{multiset1}
\displaystyle\bigcup_{\q\in \mathcal{M}^{'}(\p)}
\displaystyle\bigcup_{\p_1\in W^{'}(\q)}
T_{\p_1}\setminus\{T_1\},
\end{equation}
and the tangent line $\ell$ other than $T_1$ 
through an external point $\p_1\not=\p$
(respectively, $\p_1=\p$) on $T_1$ occurs an 
odd (respectively, even) number of times
in (\ref{multiset1}) by Lemma~\ref{set22}, 
we obtain
\begin{equation}\label{temp2}
\begin{array}{lllll}
\displaystyle\sum_{\q\in \mathcal{M}^{'}(\p)}\displaystyle\sum_{\ell\in T(\q,\ell(\q))}\chi_{\ell}
& \equiv & \displaystyle\sum_{\q\in \mathcal{M}^{'}(\p)}\displaystyle\sum_{\p_1\in W^{'}(\q)}
\displaystyle\sum_{\ell\in T_{\p_1}\setminus\{T_1\}}\chi_{\ell}\\
{} &=& \displaystyle\sum_{\p_1\in E_{T_1}}\displaystyle\sum_{\ell\in T_{\p_1}\setminus\{T_1\}}
b_{\ell}\chi_{\ell}\\
{}&=&b_{T_2}\chi_{T_2}+\displaystyle\sum_{\ell\in T\setminus\{T_1,T_2\}}b_{\ell}\chi_{\ell}\\
{}&\equiv & \displaystyle\sum_{\ell\in T\setminus\{T_1,T_2\}}\chi_{\ell} \pmod 2,

\end{array}
\end{equation}
where $b_{\ell}$ for $\ell\in T\setminus \{T_1, T_2\}$ 
are all odd integers and $b_{T_2}$ is an even integer.

Using (\ref{temp3}), (\ref{temp2}), and the fact that
$\sum_{\ell\in T}\chi_{\ell} ={\bf 0} \pmod 2$,
we have
\begin{equation*}
\begin{array}{llll}
\chi_{T_1}+\chi_{T_2} & \equiv & \displaystyle\sum_{\chi\in T\setminus\{T_1, T_2\}}\chi_{\ell}\\
{} & \equiv & \displaystyle\sum_{\q\in\mathcal{}M^{'}(\p)}
\displaystyle\sum_{\ell\in Pa_{\q}}\chi_{\ell}\\
{} & \equiv & \displaystyle\sum_{\q\in \mathcal{M}^{'}(\p)}\chi_{N_{Pa, E}(\q)} \pmod 2.
\end{array}
\end{equation*}
\QED

\section{The Conjugacy Classes and Intersection Parity}
In this section, we review the conjugacy classes of
$H$ and study their intersections with some special 
subsets of $H$.
\subsection{Conjugacy classes}
Recall that \begin{displaymath} H=
\left.\left\{\left(\begin{array}{ccc}
a^2 & ab & b^2 \\
2ac & ad+bc & 2bd \\
c^2 & cd & d^2 \\
\end{array}\right)\right| a, b, c, d\in \Ff_q,\;ad-bc=1\right\}
\end{displaymath}
is the subgroup of $G$ that is isomorphic to 
$\PSL(2,q)$. If we define $T=\text{tr}(g)+1$, 
where $g\in H$ and tr$(g)$ is the trace of $g$,
then the conjugacy classes of $H$ can be read 
as follows.

\begin{Lemma}\cite[Lemma 3.2]{swx}\label{classes}
The conjugacy classes of $H$ are given as follows.
\begin{itemize}
\item[(i)] $D=\{{\bf d}(1,1,1)\}$;
\item[(ii)] $F^{+}$ and $F^{-}$, where $F^{+}\cup F^{-} = 
\{g\in H\mid T(g) = 4,\;g\not={\bf d}(1,1,1)\}$;
\item[(iii)] $[\theta_i] =\{g\in H \mid T(g) = \theta_i\}$, 
$1\le i \le \frac{q-5}{4}$ if $q\equiv 1\pmod 4$, or 
$1\le i\le \frac{q-3}{4}$ if $q\equiv 3 \pmod 4$, where
$\theta_i\in \Sq_q$, $\theta_i\not= 4$, and 
$\theta_i-4\in\Sq_q$;
\item[(iv)] $[0]=\{g\in H\mid T(g) =0\}$;
\item[(v)] $[\pi_k]=\{g\in H\mid T(g) = \pi_k\}$, 
$1\le k\le \frac{q-1}{4}$ if $q\equiv 1\pmod 4$, 
or $1\le k\le \frac{q-3}{4}$ if $q\equiv 3 \pmod 4$, 
where $\pi_i\in\Sq_q$, $\pi_k\not=4$, and 
$\pi_k-4\in\Nsq_q$.
\end{itemize}
\end{Lemma}
\begin{Remark}\label{order}
The set $F^+\cup F^-$
forms one conjugacy class of $G$, and splits 
into two equal-sized classes $F^{+}$ and $F^{-}$ 
of $H$. For our purpose, we denote $F^+\cup F^-$ 
by $[4]$. Also, each of $D$, $[\theta_i]$, $[0]$, 
and $[\pi_k]$ forms a single conjugacy class of 
$G$. The class $[0]$ consists of all the elements 
of order $2$ in $H$.
\end{Remark}
In the following, for convenience, we frequently 
use $C$ to denote any one of $D$, $[0]$, $[4]$, 
$[\theta_i]$, or $[\pi_k]$. That is,
\begin{equation}\label{defC}
C=D, [0], [4], [\theta_i], \;{\rm or}\;[\pi_k].
\end{equation}

\subsection{Intersection properties}

\begin{Definition}
Let $\p \in I$, $\q \in E$, $\ell \in Pa$. We 
define
\begin{equation*}
\begin{array}{ccc}
\mathcal{H}_{\p,\q}=\{h\in H\mid (\p^\perp)^h\in Pa_{\q}\}
\;\;\textrm{and}\;\;
\mathcal{S}_{\p,\ell}=\{h\in H\mid (\p^\perp)^h=\ell\}.
\end{array}
\end{equation*}
That is, $\mathcal{H}_{\p,\q}$ consists of all 
the elements of $H$ that map the passant line 
$\p^\perp$ to a passant line through $\q$ and 
$\mathcal{S}_{\p,\ell}$ is the set of elements 
in $H$ that map $\p^\perp$ to the passant line 
$\ell$.
\end{Definition}
Using the above notation, since $G$ preserves
incidence, for $g\in G$, $\p\in I$, and $\ell\in Pa$,
we have
\begin{equation}\label{interest}
\mathcal{H}_{\p,\q}^g=\mathcal{H}_{\p^g,\q^g},\;
\mathcal{S}_{\p,\ell}^g=\mathcal{S}_{\p^g,\ell^g}.
\end{equation}
The following corollary is apparent.

\begin{Corollary}\label{y1}
Let $g\in G$ and $C$ be given in $($\ref{defC}$)$ 
and let $\p$ and $\q$ be two external points. Then 
$(C\cap \mathcal{H}_{\p,\q})^g=C\cap \mathcal{H}_{\p^g, \q^g}$.

\end{Corollary}

Next the size of the intersection of each conjugacy 
class of $H$ with $K$ which stabilizes an element of 
$I$ in $H$ is calculated.
\begin{Corollary}\label{y11}
Let $\p \in I$ and $K=H_\p$. Then we have
\begin{itemize}
\item[(i)] $|K\cap D|=1$;
\item[(ii)] $|K\cap [4]| = 0$;
\item[(iii)] $|K\cap [\pi_k]| = 2$ for each $k$;
\item[(iv)] $|K\cap [\theta_i]| = 0$ for each $i$;
\item[(v)] $|K\cap [0]|=\frac{q+1}{2}$ or $\frac{q-1}{2}$ 
accordingly as $q\equiv 1\pmod 4$ or $q\equiv 3\pmod 4$.
\end{itemize}
\end{Corollary}
{\Proof} Let $\q=(1,0,-\xi)$ and $K_1=H_\q$. Since 
$H$ is transitive on $I$, it follows $\q^g=\p$ for 
some $g\in H$. By Lemma~\ref{stabilizer}, we have
$K_1^g=K$. Consequently, $$|K\cap C|=|(K_1\cap C)^g|.$$
Therefore, to prove the corollary, it is enough to 
consider $\p=\q$.
It is clear that $|D\cap K|=1$. Let $g\in K\cap C$. 
Then the quadruples $(a,b,c,d)$ determining $g$ 
satisfy the following equations
\begin{equation}\label{linear1}
\begin{array}{ccccc}
bd-ac\xi& = & 0\\
b^2-a^2\xi& = & -\xi(d^2-c^2\xi)\\
ad-bc& = & 1\\
a+d & = & s,
\end{array}
\end{equation}
where $s^2=0$, $4$, $\pi_k$, $\theta_i$. The equations 
in (\ref{linear1}) give (1) $a=d=\frac{s}{2}$, 
$c^2=\frac{s^2-4}{4\xi}$, $b^2=\frac{(s^2-4)\xi}{4}$
and (2) $a=-d$, $s=0$, $c^2\xi-1=a^2$. From Case (1), 
we see that $|K\cap [\pi_k]|=2$ for each $[\pi_k]$ and 
$|K\cap C|=0$ for $C=[\theta_i]$, $[4]$; moreover, if 
$q\equiv 3\pmod 4$, we obtain one group element 
{\bf ad}$(-\xi, -1,\xi^{-1})\in K\cap [0]$ in Case (1). 
Since the number of $t\in\Nsq_q$ satisfying $t-1\in \Sq_q$ 
is $\frac{q-1}{4}$ or $\frac{q-3}{4}$ accordingly as 
$q\equiv 1\pmod 4$ or $q\equiv 3\pmod 4$ by 
Lemma~\ref{bsize}, the number of $c\in\Ff_q^*$ satisfying 
$c^2\xi-1\in\Sq_q$ is $2|(\Nsq_q-1)\cap\Sq_q|$ which is 
$\frac{q-1}{2}$ or $\frac{q-3}{2}$ accordingly as 
$q\equiv 1\pmod 4$ or $q\equiv 3\pmod 4$. When 
$q\equiv 1\pmod 4$, $c=0$ also satisfies $c^2\xi-1\in\Sq_q$.
Therefore, Case (1) and Case (2) give $\frac{q+1}{2}$ 
or $\frac{q-1}{2}$ different group elements in $K\cap [0]$ 
depending on $q$. Now the corollary is proved.
\QED

In the following lemma, we investigate the parity 
of $|\mathcal{H}_{\p,\q}\cap C|$ for each $C\not=[0]$ 
and $\p\in I$, $\q\in E$. Recall that $\ell_{\p,\q}$ 
is the line through $\p$ and $\q$.
\begin{Lemma}\label{m1}
Assume that $q\equiv 1\pmod 4$. Let $\p\in I$ and 
$\q\in E$. Suppose that $C=D$, $[4]$, $[\pi_k]$ 
$(1\le k\le \frac{q-1}{4})$, $[\theta_i]$ 
$(1\le i\le \frac{q-5}{4})$.
\begin{itemize}
\item[(i)] If $\ell_{\p,\q}\in Se_{\p}$, then 
$|\mathcal{H}_{\p,\q}\cap C|$ is even for each $C$.

\item[(ii)] If $\ell_{\p,\q}\in Pa_\p$ and 
$\q\in \p^\perp$, then $|\mathcal{H}_{\p,\q}\cap C|$ 
is odd if and only if $C=D$.

\item[(iii)] If $\ell_{\p,\q}\in Pa_{\p}$ and 
$\q\notin\p^\perp$, for each class $[\pi_k]$
with $1\le k\le \frac{q-1}{4}$, there are two 
different points $\q_1$, $\q_2\in E_{\ell_{\p,\q}}$
such that $|[\pi_k]\cap \mathcal{H}_{\p,\q_j}|$ is 
odd for $j=1$, $2$; moreover, the two points 
associated with one class $[\pi_{k_1}]$ are different 
from those associated with the other class $[\pi_{k_2}]$, 
where $[\pi_{k_1}]\not=[\pi_{k_2}]$.

\end{itemize}
\end{Lemma}
{\Proof} Since $G$ acts transitively on $I$ and 
preserves incidence, without loss of generality, 
we may assume that $\p=(1,0,-\xi)$. From (\ref{groupg}), 
it follow that 
\begin{equation}\label{stabP}
\begin{array}{lllllll}
K:=G_\p=\\
\left.\left\{\left(\begin{array}{ccc} d^2 & cd\xi & c^2 \xi^2\\ 
2cd & d^2+c^2\xi & 2cd\xi\\ c^2 & cd & d^2 \end{array}\right)\right| 
d,c \in \Ff_q, d^2-c^2\xi =1\right\}\\
\bigcup\left.\left\{\left(\begin{array}{ccc} d^2 & -cd\xi & c^2 \xi^2\\ 
2cd& -d^2-c^2\xi& 2cd\xi\\ c^2 & -cd & d^2 \end{array}\right)\right| 
d,c \in \Ff_q, -d^2+c^2\xi =1\right\}\\
\bigcup\left.\left\{\left(\begin{array}{ccc} d^2 & cd & c^2\\ 
2cd\xi^{-1} & d^2+c^2\xi^{-1}& 2cd\\ c^2\xi^{-2} & cd\xi^{-1} & d^2 
\end{array}\right)\right| d,c \in \Ff_q, d^2\xi-c^2 =1\right\}\\
\bigcup\left.\left\{\left(\begin{array}{ccc} d^2 & -cd & c^2 \\ 
2cd\xi^{-1} & -d^2-c^2\xi^{-1}& 2cd\\ c^2\xi^{-2} & -cd\xi^{-1} & d^2 
\end{array}\right)\right| d,c \in \Ff_q, -d^2\xi+c^2 =1\right\}.\\
\end{array}
\end{equation}
Since $K$ is transitive on both $Pa_{\p}$ and $Se_{\p}$ 
by Proposition~\ref{Ktransitive} and 
$$|\mathcal{H}_{\p,\q}\cap C| = 
|(\mathcal{H}_{\p,\q}\cap C)^g| = |\mathcal{H}_{\p,\q^g}\cap C|$$
by Corollary~\ref{y1}, we may assume that $\q$ is on 
either $\ell_1$ or $\ell_2$, where
$\ell_1=[1,0,\xi^{-1}]\in Pa_{\p}$ and 
$\ell_2=[0,1,0]\in Se_{\p}$.

{\bf Case I.} $\q\in \ell_1$ and 
$\q\notin\p^\perp$.

In this case, $\q=(1,x,-\xi)$ for some $x\in\Ff_q^*$ 
and $x^2+\xi\in\Sq_q$, and 
$$Pa_{\q}=\{[1,s,(1+sx)\xi^{-1}]\mid s\in \Ff_q, 
s^2-4(1+sx)\xi^{-1}\in\Nsq_q\}.$$
Using (\ref{stabP}), we obtain that
$$K_{\q}=\{{\bf d}(1,1,1), {\bf ad}(1,-\xi^{-1},\xi^{-2})\}.$$
It is apparent that ${\bf d}(1,1,1)$ fixes each line 
in $Pa_{\q}$. From 
$${\bf ad}(1,-\xi^{-1},\xi^{-2})^{-1}(1,s,(1+sx)\xi^{-1})^\top=((1+sx)\xi, -s\xi, 1)^\top,$$
it follows that $[1,s,(1+sx)\xi^{-1}]\in Pa_{\q}$ is 
fixed by $K_{\q}$ if and only if $s=0$ or $s=-2x^{-1}$. 
Therefore, $K_{\q}$ has two orbits of length $1$ on 
$Pa_\q$,
i.e. 
$\{\ell_1=[1,0,\xi^{-1}]\}$ and 
$\{\ell_3=[1,-2x^{-1}, -\xi^{-1}]\}$, 
and all other orbits, whose representatives are 
$\mathcal{R}_1$, have length $2$. From
\begin{equation*}
|\mathcal{H}_{\p,\q}\cap C|=|\mathcal{S}_{\p,\ell_1}\cap C|
+|\mathcal{S}_{\p,\ell_3}\cap C|+ 
2\displaystyle\sum_{\ell\in\mathcal{R}_1}
|\mathcal{S}_{\p,\ell}\cap C|,
\end{equation*}
it follows that the parity of $|\mathcal{H}_{\p,\q}\cap C|$ 
is determined by that of $|\mathcal{S}_{\p,\ell_1}\cap C|+
|\mathcal{S}_{\p,\ell_3}\cap C|$. Here we used the fact that 
$|\mathcal{S}_{\p,\ell}\cap C|=|\mathcal{S}_{\p,\ell^{'}}\cap C|$ 
if $\{\ell, \ell^{'}\}$ is an orbit of $K_{\p}$ on $Pa_{\q}$. 
It is clear that
$|\mathcal{S}_{\p,\ell}\cap D|=|\mathcal{S}_{\p,\ell_3}\cap D|=0$.

Note that the quadruples $(a,b,c,d)$ that determine group 
elements in $\mathcal{S}_{\p,\ell_1}\cap C$ satisfy the 
following equations
\begin{equation}\label{s1}
\begin{array}{ccccc}
-2cd+2ab\xi^{-1} & = & 0\\
c^2-a^2\xi^{-1} & = & (d^2-b^2\xi^{-1})\xi^{-1}\\
a+d & = & s\\
ad-bc & = & 1\\

\end{array}
\end{equation}
where $s^2=4$, $\pi_k$, $\theta_i$. The first two 
equations in (\ref{s1})give $c=\pm\sqrt{-1}c\xi^{-1}$
and $a=\pm\sqrt{-1} d$. Combining them with the 
last two equationsin (\ref{s1}), we obtain $0$, $4$ 
or $8$ quadruples $(a,b,c,d)$ satisfying the above 
equations, among which, both $(a,b,c,d)$ and 
$(-a,-b,-c,-d)$ do appear at the same time. Therefore, 
$|\mathcal{S}_{\p,\ell_1}\cap C|$ is $0$, $2$, or $4$.
Particularly, in $[0]$, there might be only $2$ 
elements satisfying the above conditions.

Similarly, the quadruples $(a,b,c,d)$ that determine 
a group element in $\mathcal{S}_{\p,\ell_3}\cap C$ 
satisfy the following equations
\begin{equation}\label{s2}
\begin{array}{cccccccc}
-2cd+2ab\xi^{-1} & = & -2x^{-1}(d^2-b^2\xi^{-1})\\
c^2-a^2\xi^{-1} & = & -\xi^{-1}(d^2-b^2\xi^{-1})\\
a+d & = & s\\
ad-bc & =& 1,
\end{array}
\end{equation}
where $s^2=4$, $\pi_k$, $\theta_i$. The first two 
equations in (\ref{s2}) give
\begin{equation}\label{s3}
d^2-b^2\xi^{-1}=\pm A,
\end{equation}
where 
\begin{equation}\label{s4}
A=\sqrt{\frac{1}{(x^2+\xi^{-1})\xi}}.
\end{equation}
From (\ref{s3}), $c^2-d^2=\mp A\xi^{-1}$ and 
$a^2+d^2=s^2-2-2bc$, it follows
that
\begin{equation}\label{s5}
(b\xi^{-1}+c)^2=-(\pm2A+2-s^2)\xi^{-1}.
\end{equation}
Hence, if (\ref{s2}) determines an odd number 
of group elements, then
\begin{equation*}
-(\pm 2A+2-s^2)\xi^{-1}\notin\Nsq_q.
\end{equation*}

If $-(\pm2A+2-s^2)\xi^{-1}\in \Sq_q$ and we set 
$B_{(\pm)}:=\sqrt{-(\pm 2A+2-s^2)\xi^{-1}}$, by 
$c^2-d^2\xi^{-1}=\pm A\xi^{-1}$ and 
$a^2-b^2\xi^{-1}=\mp A\xi^{-1}$, we have
\begin{equation}\label{s6}
d=\frac{1}{2s}[s^2+(\xi B_{(\pm)}^2-2B_{(\pm)}b)]\;\left(\text{or}\;
d=\frac{1}{2s}[s^2+(\xi B_{(\pm)}^2+2B_{(\pm)}b)]\right)
\end{equation}
and thus
\begin{equation}\label{s7}
a=\frac{1}{2s}[s^2-(\xi B_{(\pm)}^2-2B_{(\pm)}b)]\;\left(\text{or}\;
a=\frac{1}{2s}[s^2-(\xi B_{(\pm)}^2+2B_{(\pm)}b)]\right).
\end{equation}
Combining $b=(\pm B_{(\pm)}-c)\xi^{-1}$ and $ad-bc=1$, we have
\begin{equation}\label{s8_1}
\left(\xi - \frac{B_{(\pm)}^2\xi^2}{s^2}\right)c^2+
\left(\frac{\xi B_{(\pm)}^3\xi^2}{s^2}-B_{(\pm)}\xi\right)c +
\left(\frac{s^2}{4}-\frac{B_{(\pm)}^4\xi^2}{4s^2}-1\right)=0
\end{equation}
or
\begin{equation}\label{s8_2}
\left(\xi - \frac{B_{(\pm)}^2\xi^2}{s^2}\right)c^2-
\left(\frac{\xi B_{(\pm)}^3\xi^2}{s^2}-B_{(\pm)}\xi\right)c +
\left(\frac{s^2}{4}-\frac{B_{(\pm)}^4\xi^2}{4s^2}-1\right)=0.
\end{equation}
The discriminant of (\ref{s8_1}) or (\ref{s8_2}) is
\begin{equation}\label{s9}
\Delta=\xi^2\left(\frac{\xi B_{(\pm)}^2\xi}{s^2}-1\right)
\left(\frac{s^2}{\xi}-\frac{4}{\xi}-B_{(\pm)}^2\right)
=\frac{4x^2\xi}{(x^{-2}+\xi^{-1})s^2}\in\Sq_q.
\end{equation}
From (\ref{s6}), (\ref{s7}), (\ref{s8_1}), (\ref{s8_2}),
and (\ref{s9}), it follows that (\ref{s2}) produces $2$ or 
$4$ group elements; that is,
$|\mathcal{S}_{\p,\ell_3}\cap C|=2$ or $4$.

If $-(\pm 2A+2-s^2)\xi^{-1}=0$, then $s^2$ is one of 
$2A+2$ and $-2A+2$ since
\begin{equation*}
(2A+2)(-2A+2)=\frac{4x^2}{x^2+\xi^{-1}}\in \Nsq_q.
\end{equation*}
Therefore, in this case, we have 
$|[s^2]\cap \mathcal{S}_{\p,\ell_3}|=1$. It is also clear 
that, for the same $[s^2]$, 
$$|[s^2]\cap\mathcal{S}_{\p,\p_1^\perp}|$$
is odd, where $\p_1=(1,-x,-\xi)\in E_{\ell_3}$. Moreover, 
when $x$ runs over $$L:=\{x\in\Ff_q^*\mid x^2+\xi\in\Sq_q\}$$ 
once, each $[\pi_k]$ with $1\le k\le \frac{q-1}{4}$ appears 
exactly twice in the multiset
\begin{equation*}
\left\{2\sqrt{\frac{1}{(x^{-2}+\xi^{-1})\xi}}+2\right\}
\bigcup
\left\{-2\sqrt{\frac{1}{(x^{-2}+\xi^{-1})\xi}}+2\right\}.
\end{equation*}
Note that 
\begin{equation*}
\pm \frac{2}{\sqrt{(x_1^{-2}+\xi^{-1})\xi}}+2=
\pm \frac{2}{\sqrt{(x_2^{-2}+\xi^{-1})\xi}}+2
\end{equation*}
if and only if $x_1=\pm x_2$. Therefore, for each class 
$[\pi_k]$ with $1\le k\le \frac{q-1}{4}$,
there are two different points $\q_1$, $\q_2\in E_{\ell_{\p,\q}}$ 
such that $|[\pi_k]\cap\mathcal{H}_{\p,\q_j}|$ is odd for 
$j=1$, $2$; further, the two points associated with one 
class $[\pi_{k_1}]$ are different from those associated 
with the other class $[\pi_{k_2}]$, where 
$[\pi_{k_1}]\not=[\pi_{k_2}]$. The proof of (iii) is 
completed.

{\bf Case II.} $\q=\ell_1\cap \p^\perp$.

In this case, $\q=(0,1,0)$. From (\ref{stabP}), it follows that
\begin{equation*}
K_{\q}=\{{\bf d}(1,1,1), {\bf ad}(-1,1,-1),
{\bf d}(-1,-\xi^{-1},-\xi^{-2}), 
{\bf ad}(1,-\xi^{-1},\xi^{-2})\}.
\end{equation*}
Since $Pa_{\q}=\{[1,0,-x]\mid x\in\Nsq_q\}$, it follows that 
the passant lines through $\q$ that are fixed by $K_\q$ are 
$\ell_1=[1,0,\xi^{-1}]$ and $\ell_4=[1,0,-\xi^{-1}]$. Thus, 
$K_{\q}$ has two orbits of length $1$ on $Pa_{\q}$ and all 
the other orbits, whose representatives are $\mathcal{R}_2$, 
have length $2$. By
\begin{equation*}
|\mathcal{H}_{\p,\q}\cap C|=|\mathcal{S}_{\p,\ell_1}\cap C|
+|\mathcal{S}_{\p,\ell_4}\cap C|+2\displaystyle
\sum_{\ell\in\mathcal{R}_2}|\mathcal{S}_{\p,\ell}\cap C|,
\end{equation*}
we obtain that the parity of $|\mathcal{H}_{\p,\q}\cap C|$ 
is determined by that of $|\mathcal{S}_{\p,\ell_1}\cap C|+
|\mathcal{S}_{\p,\ell_4}\cap C|$. From the discussions
in {\bf Case I}, we know that $|\mathcal{S}_{\p,\ell_1}\cap C|$ 
always even. Since $\ell_4=\p^\perp$ and $G_{\p}=G_{\p^\perp}$ 
by Lemma~\ref{a11}, it follows from Corollary~\ref{y11} 
that $|\mathcal{S}_{\p,\ell_4}\cap C|$ is odd if and only 
if $C=D$. The proof of (ii) is completed.

{\bf Case III.} $\q\in\ell_2$.

In this case , $\q=(1,0,-y)$ for some $y\in\Sq_q$. 
Using (\ref{stabP}), we see that
$$K_{\q}=\{{\bf d}(1,1,1), {\bf d}(-1,1,-1)\}.$$
Moreover, all the orbits of $K_{\q}$ on 
$Pa_{\q}=\{[1,s,y^{-1}]\mid s\in\Ff_q^*, s^2-4y^{-1}\in\Nsq_q\}$
have length $2$, then $|\mathcal{H}_{\p,\q}\cap C|$ 
is even for each $C$. Part (i) is proved.

\QED

\begin{Lemma}\label{m2}
Assume that $q\equiv 3\pmod 4$. Let $\p\in I$ and 
$\q \in E$. Suppose that $C=D$, $[4]$, $[\pi_k]$ 
$(1\le k\le \frac{q-3}{4})$, $[\theta_i]$ 
$(1\le i\le \frac{q-3}{4})$.

\begin{itemize}
\item[(i)] If $\ell_{\p,\q}\in Se_{\p}$ and $\q\notin\p^\perp$, 
for each class $[\theta_i]$ with $1\le i\le \frac{q-3}{4}$, 
there are two different points $\q_1$, $\q_2\in E_{\ell_{\p,\q}}$ 
such that $|[\theta_i]\cap \mathcal{H}_{\p,\q_j}|$ is odd for $j=1$, 
$2$; moreover, the two points associated with one class 
$[\theta_{i_1}]$ are different from those associated with the 
other class $[\theta_{i_2}]$, where 
$[\theta_{i_1}]\not=[\theta_{i_2}]$.
\item[(ii)] If $\ell_{\p,\q}\in Se_\p$ and $\q\in \p^\perp$, 
then $|\mathcal{H}_{\p,\q}\cap C|$ is odd if and only if $C=D$.
\item[(iii)] If $\ell_{\p,\q}\in Pa_{\p}$, then 
$|\mathcal{H}_{\p,\q}\cap C|$ is even for each $C$.

\end{itemize}

\end{Lemma}
{\Proof} The proof is essentially the same as the one of 
Lemma~\ref{m1}. 
We omit the details.\QED

\section{Group Algebra $FH$}
\subsection{2-Blocks of H}
In this section we recall several results on the $2$-blocks 
of $H\cong PSL(2,q)$. We refer the reader to {\rm \cite{gabriel}} 
or {\rm{\cite{brauer}}} for a general introduction on 
this subject.

Let $\mathbf{R}$ be the ring of algebraic integers in 
the complex field $\Cc$. We choose a maximal ideal 
$\mathbf{M}$ of $\mathbf{R}$ containing $2\mathbf{R}$. 
Let $F=\mathbf{R}/\mathbf{M}$ be the residue field of 
characteristic $2$, and let $* :\mathbf{R}\rightarrow F$ 
be the natural ring homomorphism. Define
\begin{equation}
\begin{array}{llll}\label{ring_s}
\mathbf{S}& =& \{\frac{r}{s}\mid r\in\mathbf{R},\;s\in\mathbf{R}\setminus\mathbf{M}\}.\\
\end{array}
\end{equation}
Then it is clear that the map $* : \mathbf{S}\rightarrow F$ defined
by $(\frac{r}{s})^* = r^*(s^*)^{-1}$ is a ring homomorphism with
kernel $\mathcal{P} = \{\frac{r}{s}\mid
r\in\mathbf{M},\;s\in\mathbf{R}\setminus\mathbf{M}\}$. In the rest
of this article, $F$ will always be the field of characteristic $2$
constructed as above. Note that $F$ is an algebraic closure of
$\Ff_2$.

Let $Irr(H)$ and $IBr(H)$ be the set of irreducible ordinary
characters and the set of irreducible Brauer characters of $H$,
respectively. In the following, we deduce the $2$-blocks of $H$
from the known results on the $2$-blocks of $\PSL(2,q)$. For
baisc results on blocks of finite groups, we refer the reader to
Chapter 3 of \cite{gabriel}.

The character tables of $\PSL(2,q)$ were obtained by Jordan and 
Schur independently; see\cite{jordan}, \cite{GK}, or \cite{schur} 
for the detailed discussions. The irreducible characters of $H$ 
can be read off from the character tables of $\PSL(2,q)$ as follows.

\begin{Lemma}$($\cite{jordan}, \cite{GK}, \cite{schur}$)$
The irreducible ordinary characters of $H$ are:
\begin{itemize}
\item[(i)] $1=\chi_0$, $\gamma$, $\chi_1$, ..., $\chi_{\frac{q-1}{4}}$,
$\beta_1$, $\beta_2$, $\phi_1$, ..., $\phi_{\frac{q-5}{4}}$ if 
$q\equiv 1\pmod 4$, where $1=\chi_0$ is the trivial character,
$\gamma$ is the character of degree $q$, $\chi_s$ for $1\le s\le \frac{q-1}{4}$
are the characters of degree $q-1$, $\phi_r$ for $1\le r\le \frac{q-5}{4}$
are the characters of degree $q+1$, and $\beta_i$ for $i=1$, $2$ are the
characters of degree $\frac{q+1}{2}$;
\item[(ii)]  $1=\chi_0$, $\chi_1$, ..., $\chi_{\frac{q-3}{4}}$,
$\beta_1$, $\eta_2$, $\eta_1$, ..., $\phi_{\frac{q-3}{4}}$ if 
$q\equiv 3\pmod 4$, where $1=\chi_0$ is the trivial character,
$\gamma$ is the character of degree $q$, $\chi_s$ for $1\le s\le \frac{q-1}{3}$
are the characters of degree $q-1$, $\phi_r$ for $1\le r\le \frac{q-3}{4}$
are the characters of degree $q+1$, and $\eta_i$ for $i=1$, $2$ are the
characters of degree $\frac{q-1}{2}$;
\end{itemize}
\end{Lemma}
The following lemma tells us how the irreducible ordinary 
characters of $H$ are partitioned into $2$-blocks.
\begin{Lemma}\cite[Lemma 4.1]{swx}\label{blocks}
First assume that $q\equiv 1 \pmod 4$ and $q-1 = m2^n$, 
where $2\nmid m$.
\begin{itemize}

\item[(i)] The principal block $B_0$ of $H$ contains 
$2^{n-2}+3$ irreducible characters
$$\chi_0 = 1,\;\gamma,\;\beta_1,\;\beta_2,\;\phi_{i_1}, ..., 
\phi_{i_{(2^{n-2}-1)}},$$ where $\chi_0=1$ is the trivial 
character of $H$, $\gamma$ is the irreducible character of 
degree $q$ of $H$,  $\beta_1$ and $\beta_2$ are the two 
irreducible characters of degree $\frac{q+1}{2}$, and 
$\phi_{i_k}$ for $1\le k\le 2^{n-2}-1$ are distinct 
irreducible characters of degree $q+1$ of $H$.

\item[(ii)]  $H$ has $\frac{q-1}{4}$ blocks $B_s$ of 
defect $0$ for $1\le s \le \frac{q-1}{4}$, each of which 
contains an irreducible ordinary character $\chi_s$ of 
degree $q-1$.

\item[(iii)] If $m\ge 3$, then $H$ has $\frac{m-1}{2}$ 
blocks $B_t^{'}$ of defect $n-1$ for $1\le t\le \frac{m-1}{2}$, 
each of which contains $2^{n-1}$ irreducible ordinary 
characters $\phi_{t_i}$ for $1\le i \le 2^{n-1}$.

\end{itemize}

Now assume that $q\equiv 3\pmod 4$ and $q+1=m2^n$, 
where $2\nmid m$ .
\begin{itemize}

\item[(iv)] The principal block $B_0$ of $H$ contains 
$2^{n-2}+3$ irreducible characters
$$\chi_0 = 1,\;\gamma,\;\eta_1,\;\eta_2,\;\chi_{i_1}, ..., \chi_{i_{(2^{n-2}-1)}},$$
where $\chi_0=1$ is the trivial character of $H$, $\gamma$ 
is the irreducible character of degree $q$ of $H$, $\eta_1$ 
and $\eta_2$ are the two irreducible characters of degree
$\frac{q-1}{2}$, and $\chi_{i_k}$ for $1\le k\le 2^{n-2}-1$ 
are distinct irreducible characters of degree $q-1$ of $H$.

\item[(v)] $H$ has $\frac{q-3}{4}$ blocks $B_r$ of defect 
$0$ for $1\le r\le \frac{q-3}{4}$, each of which contains 
an irreducible ordinary character $\phi_r$ of degree $q+1$.

\item[(vi)] If $m\ge 3$, then $H$ has $\frac{m-1}{2}$ blocks 
$B_t^{'}$ of defect $n-1$ for $1\le t\le \frac{m-1}{2}$, 
each of which contains $2^{n-1}$ irreducible ordinary 
characters $\chi_{t_i}$ for $1\le i \le 2^{n-1}$.

\end{itemize}

Moreover, the above blocks form all the $2$-blocks of $H$.
\end{Lemma}
\begin{Remark}
Parts $($i$)$ and $($iv$)$ are from Theorem 1.3 in 
{\rm\cite{landrock}} and their proofs can be found in 
Chapter 7 of {III} in {\rm \cite{brauer}}. Parts $($ii$)$ 
and $($v$)$ are special cases of Theorem 3.18 in 
{\rm \cite{gabriel}}. Parts $($iii$)$ and $($vi$)$ are
proved in Sections {II} and {VIII} of {\rm \cite{burkhardt}}.
\end{Remark}


\subsection{Block Idempotents}

Let $Bl(H)$ be the set of $2$-blocks of $H$. If $B\in Bl(H)$, 
we write $$f_B = \displaystyle\sum_{\chi\in Irr(B)}e_{\chi},$$
where $e_{\chi}=\frac{\chi(1)}{|H|}\sum_{g\in H} \chi(g^{-1})g$ 
is a central primitive idempotent of $\mathbf{Z}(\Cc H)$ and
$Irr(B)=Irr(H)\cap B$. For future use, we define $IBr(B)=IBr(H)
\cap B$. Since $f_B$ is an element of $\mathbf{Z}(\Cc H)$, 
we may write
\begin{displaymath}
\begin{array}{lllll}
f_B& =& \displaystyle\sum_{C\in cl(H)}f_B(\widehat{C})\widehat{C},\\
\end{array}
\end{displaymath}
where $cl(H)$ is the set of conjugacy classes of $H$, 
$\widehat{C}$ is the sum of elements in the class $C$, and
\begin{equation}\label{id}
\begin{array}{llll}
f_B(\widehat{C})& = &\frac{1}{|H|}\displaystyle\sum_{\chi\in Irr(B)}
\chi(1)\chi(x_C^{-1})
\end{array}
\end{equation}
with a fixed element $x_C\in C$.
\begin{Theorem}\label{osima}
Let $B\in Bl(H)$. Then $f_B\in \mathbf{Z}(\mathbf{S}H)$. 
In other words,$f_B(\widehat{C})\in \mathbf{S}$ for each 
block of $H$.

\end{Theorem}
{\Proof} It follows from Corollary 3.8 in {\rm \cite{gabriel}}. 
\QED

We extend the ring homomorphism $*: \mathbf{S}\rightarrow F$ 
to a ring homomorphism $*:\mathbf{S}H\rightarrow FH$ by setting
$(\sum_{g\in H} s_g g)^*= \sum_{g\in H} s_g^* g$. Note that $*$ 
maps $\mathbf{Z}(\mathbf{S}H)$ onto $\mathbf{Z}(FH)$ via 
$(\sum_{C\in cl(H)}s_C \widehat{C})^*$ = $\sum_{C\in cl(H)} 
s_C^* \widehat{C}$. Now we define
$$e_B = (f_B)^* \in \mathbf{Z}(FH),$$
which is the {\it block idempotent} of $B$. Note that $e_B
e_{B^{'}}= \delta_{B B^{'}}e_B$ for $B$, $B^{'}\in Bl(H)$, 
where $\delta_{B B^{'}}$ equals 1 if $B=B'$, 0 otherwise. 
Also $1=\sum_{B\in Bl(H)}e_B$.

All the block idempotents of the $2$-blocks of $H$ are given 
in the following lemma; see \cite{swx} for the detailed 
calculations.
\begin{Lemma}\cite[Lemma 4.4]{swx}\label{expression}
First assume that $q\equiv 1 \pmod 4$ and $q-1= m2^n$ with 
$2\nmid m$.
\begin{itemize}

\item[1.] Let $B_0$ be the principal block of $H$. Then
\begin{enumerate}
\renewcommand{\labelenumi}{(\alph{enumi})}
\item 
$e_{B_0}(\widehat{D})=1$.

\item 
 $e_{B_0}(\widehat{F^+})= e_{B_0}(\widehat{F^-})\in F$.

\item 
$e_{B_0}(\widehat{[\theta_i]})\in F$, $e_{B_0}(\widehat{[0]})=0$.

\item 
$e_{B_0}(\widehat{[\pi_k]})=1$.
\end{enumerate}
\item[2.] Let $B_s$ be any block of defect $0$ of $H$. 
Then
\begin{enumerate}
\renewcommand{\labelenumi}{(\alph{enumi})}
\item 
$e_{B_s}(\widehat{D})=0$.
\item 
$e_{B_s}(\widehat{F^+})= e_{B_s}(\widehat{F^-})=1$.
\item 
$e_{B_s}(\widehat{[0]})= e_{B_s}(\widehat{[\theta_i]})= 0$.
\item 
$e_{B_s}(\widehat{[\pi_k]})\in F$.
\end{enumerate}
\item[3.] Suppose $m\ge 3$ and let $B_t^{'}$ be any block 
of defect $n-1$ of $H$. Then
\begin{enumerate}
\renewcommand{\labelenumi}{(\alph{enumi})}
\item 
$e_{B_t^{'}}(\widehat{D})=0$.
\item 
$e_{B_t^{'}}(\widehat{F^+})= e_{B_t^{'}}(\widehat{F^-})=1$.
\item 
$e_{B_t^{'}}(\widehat{[\theta_i]})\in F$, $e_{B_t^{'}}(\widehat{[0]})=0$.
\item 
$e_{B_t^{'}}(\widehat{[\pi_k]}) = 0$.
\end{enumerate}
\end{itemize}
Now assume that $q\equiv 3 \pmod 4$. Suppose that 
$q+1=m2^n$ with $2\nmid m$.
\begin{itemize}

\item[4.] Let $B_0$ be the principal block of $H$. 
Then
\begin{enumerate}
\renewcommand{\labelenumi}{(\alph{enumi})}
\item 
$e_{B_0}(\widehat{D})=1$.

\item 
$e_{B_0}(\widehat{F^+})= e_{B_0}(\widehat{F^-})\in F$.

\item 
$e_{B_0}(\widehat{[\theta_i]})=1$.

\item 
$e_{B_0}(\widehat{[0]})=0$, $e_{B_0}(\widehat{[\pi_k]})\in F$.
\end{enumerate}
\item[5.] Let $B_r$ be any block of defect $0$ of $H$. 
Then
\begin{enumerate}
\renewcommand{\labelenumi}{(\alph{enumi})}
\item 
$e_{B_r}(\widehat{D})=0$.
\item 
$e_{B_r}(\widehat{F^+})= e_{B_r}(\widehat{F^-})=1$.
\item 
$e_{B_r}(\widehat{[0]})= e_{B_r}(\widehat{[\pi_k]}) = 0$.
\item 
$e_{B_r}(\widehat{[\theta_i]})\in F$.
\end{enumerate}
\item[6.] Suppose that $m\ge 3$ and let $B_t^{'}$ be 
any block of defect $n-1$ of $H$. Then
\begin{enumerate}
\renewcommand{\labelenumi}{(\alph{enumi})}
\item 
$e_{B_t^{'}}(\widehat{D})=0$.
\item 
$e_{B_t^{'}}(\widehat{F^+})= e_{B_t^{'}}(\widehat{F^-})=1$.
\item 
$e_{B_t^{'}}(\widehat{[\theta_i]})= 0$.
\item 
$e_{B_t^{'}}(\widehat{[0]})=0$, $e_{B_t^{'}}(\widehat{[\pi_k]})\in F$.
\end{enumerate}

\end{itemize}
\end{Lemma}

The following corollary will be used in the proof 
of Lemma~\ref{y4}.

\begin{Corollary}\label{EX}
Let $B_s$ $($$1\le s\le \frac{q-1}{4}$$)$ or $B_r$ 
$($$1\le r\le \frac{q-3}{4}$$)$ be the blocks of 
defect $0$ of $H$ depending on whether $q\equiv 1\pmod 4$
or $q\equiv 3\pmod 4$. Using the above notation,
\begin{itemize}
\item[(i)] if $q\equiv 1 \pmod 4$, for each $B_s$, 
there is a class $[\pi_k]$ such that 
$e_{B_s}(\widehat{[\pi_k]})\not=0$;
\item[(ii)] if $q\equiv 3\pmod 4$, for each $B_r$, 
there is a class $[\theta_i]$ such that 
$e_{B_r}(\widehat{[\theta_i]})\not=0$.

\end{itemize}

\end{Corollary}
{\Proof} First we assume that $q\equiv 1\pmod 4$. 
From Theorem 8.9 in \cite{GK}, we have 
$\chi_s(g_k)=-\delta^{(2k)s}-\delta^{-(2k)s}$ for 
$1\le k\le \frac{q-1}{4}$, where $\chi_s$ is the 
irreducible ordinary character lying in $B_s$, 
$g_k\in [\pi_k]$, and $\delta$ is a primitive 
$(q+1)$-th root of unit in $\Cc$. Note that
\begin{equation*}
\begin{array}{llll}
f_{B_s}(\widehat{[\pi_k]})& = & \frac{1}{|H|}
\displaystyle\sum_{\chi_s\in B_s}
\chi_s(1)\chi_s(g_k^{-1})\\
{} & = & -\frac{q-1}{|H|}(\delta^{(2k)s}+\delta^{-(2k)s}).
\end{array}
\end{equation*}
Since 
\begin{equation*}
\begin{array}{lllll}
\displaystyle\sum_{k=1}^{(q-1)/4}e_{B_s}
(\widehat{[\pi_k]}) & = &
(-\frac{q-1}{|H|}\displaystyle\sum_{k=1}^{(q-1)/4}
\delta^{(2k)s}+\delta^{-(2k)s})^*\\
{} & = & (\frac{\delta^{2s}-\delta^{\frac{q+3}{2}s}}
{1-\delta^{2s}}+\frac{\delta^{-2s}-\delta^{-\frac{q+3}{2}s}}
{1-\delta^{-2s}})^*\\
&=& (\frac{\delta^{2s}-\delta^{\frac{q+3}{2}s}}
{1-\delta^{2s}}+\frac{\delta^{-2s}-\delta^{\frac{q-1}{2}s}}{1-\delta^{-2s}})^*\\
&=& (\frac{\delta^{2s}-\delta^{\frac{q+3}{2}s}}
{1-\delta^{2s}}+\frac{1-\delta^{\frac{q+3}{2}s}}{\delta^{2s}-1})^*\\
{} & = & 1,
\end{array}
\end{equation*}
we conclude that $e_{B_s}(\widehat{[\pi_k]})\not=0$ 
for some $k$. Part (i) is proved.

Part (ii) can be proved in the same fashion using 
Theorem 8.11 in \cite{GK}; we omit the details.\QED

Let $M$ be an $\mathbf{S}H$-module. We denote the 
reduction $M/\mathcal{P}M$, which is an $FH$-module, 
by $\overline{M}$. Then the following lemma is apparent.

\begin{Lemma}\label{reduction}
Let $M$ be an $\mathbf{S}H$-module and $B\in Bl(H)$. 
Using the above notation, we have 
$$\overline{f_B M} = e_B \overline{M},$$
i.e. reduction commutes with projection onto a block 
$B$.
\end{Lemma}

\section{Linear Maps and Their Matrices}

Let $F$ be the algebraic closure of $\Ff_2$ defined 
in Section 4. From now on, $\chi_N$ for $N\subseteq E$ 
will be always regarded as a vector over $F$. Recall 
that for $\p\in I$, $N_{Pa, E}(\p)$ 
(respectively, $N_{Se, E}(\p)$) is the set of external 
points on the passant (respectively, secant) lines 
through $\p$. We define $\mathbf{D}$ 
(respectively, $\mathbf{D}^{'}$) to be the incidence
matrix of $E$ and $N_{Pa, E}(\p)$ 
(respectively, $N_{Se, E}(\p)$)
for $\p\in I$. Namely, the columns of $\mathbf{D}$ 
and $\D^{'}$
can be viewed as the characteristic vectors of 
$N_{Pa, E}(\p)$ and $N_{Se, E}(\p)$, respectively. 
In the following, we always regard both $\D$ and 
$\D^{'}$ as matrices over $F$.

\begin{Definition}
For $\p\in I$, we define $\mathcal{G}_{\p}$ to 
be the column characteristic vector of $\p$ with 
respect to $I$, i.e. $\mathcal{G}_{\p}$ is a 
$0$-$1$ column vector of length $|I|$ with entries 
indexed by the internal points; the entry of 
$\mathcal{G}_{\p}$ is $1$ if and only if it is 
indexed by $\p$.
\end{Definition}

Let $k$ be the complex field $\Cc$, the algebraic 
closure $F$ of $\Ff_2$, or the ring  $\mathbf{S}$ 
in~(\ref{ring_s}). Let $k^I$ and $k^E$ be the free 
$k$-modules with the bases 
$\{\mathcal{G}_{\p}\mid \p\in I\}$ and 
$\{\chi_{\p}\mid \p\in E\}$, respectively. If we 
extend the actions of $H$ on the bases of $k^I$ 
and $k^E$, which are defined by 
$\chi_{\p}\cdot h=\chi_{\p^h}$ and 
$\mathcal{G}_{\q}\cdot h=\mathcal{G}_{\q^h}$
for $\p \in I$, $\q\in E$, and $h\in H$, linearly 
to $k^I$ and $k^E$ respectively, then both $k^I$ 
and $k^E$ are $kH$-permutation modules. Since $H$
is transitive on $I$, we have 
$$k^I=\textrm{Ind}_K^H(1_k),$$
where $K$ is the stabilizer of an element of $I$ 
in $H$ and $\text{Ind}_K^H(1_k)$ is the $kH$-module 
induced by $1_k$.

The decomposition of $1\uparrow_K^H$, the 
character of $\textrm{Ind}_K^H(1_k)$, into a sum 
of the irreducible ordinary characters of $H$ is 
given as follows.

\begin{Lemma}\cite[Lemma 5.2]{wu}\label{decomposition_1}
Let $K$ be the stabilizer of an internal
point in $H$. 

Assume that $q\equiv 1\pmod 4$. Let $\chi_s$,
$1\le s \le \frac{q-1}{4}$, be the irreducible 
ordinary characters of degree $q-1$, $\phi_r$, 
$1\le r \le \frac{q-5}{4}$, irreducible ordinary 
characters of degree $q+1$, $\gamma$ the 
irreducible of degree $q$, and $\beta_j$, $1\le j\le 2$,
irreducible ordinary characters of degree 
$\frac{q+1}{2}$.

\begin{itemize}
\item[(i)] If $q\equiv 1 \pmod 8$, then 
$$1\uparrow_{K}^H =1 + \displaystyle
\sum_{s=1}^{(q-1)/4}\chi_s + \gamma +\beta_1+
\beta_2 + \displaystyle\sum_{j=1}^{(q-9)/4}
\phi_{r_j},$$ where $\phi_{r_j}$, 
$1\le j\le\frac{q-9}{4}$, may not be distinct.
\item[(ii)] If $q\equiv 5\pmod 8$, then 
$$1\uparrow_{K}^H =1 + \displaystyle
\sum_{s=1}^{(q-1)/4}\chi_s + \gamma +
\displaystyle\sum_{j=1}^{(q-5)/4}\phi_{r_j},$$ 
where $\phi_{r_j}$, $1\le j\le\frac{q-5}{4}$, 
may not be distinct.
\end{itemize}

Next assume that $q\equiv 3\pmod 4$. Let 
$\chi_s$, $1\le s \le\frac{q-3}{4}$, be the 
irreducible ordinary characters of degree
$q-1$, $\phi_r$, $1\le r \le \frac{q-3}{4}$, 
the irreducible ordinary characters of degree 
$q+1$, $\gamma$ the irreducible character of 
degree $q$, and $\eta_j$, $1\le j\le 2$, the
irreducible ordinary characters of degree 
$\frac{q-1}{2}$.

\begin{itemize}
\item[(iii)] If $q\equiv 3\pmod 8$, then 
$$1\uparrow_{K}^H=
1+\displaystyle\sum_{r=1}^{(q-3)/4}\phi_r+\eta_1+
\eta_2+\displaystyle\sum_{j=1}^{(q-3)/4}\chi_{s_j},$$
where $\chi_{s_j}$, $1\le j\le\frac{q-3}{4}$, 
may not be distinct.

\item[(iv)] If $q\equiv 7\pmod 8$, then 
$$1\uparrow_{K}^H=1
+\displaystyle\sum_{r=1}^{(q-3)/4}\phi_r+
\displaystyle\sum_{j=1}^{(q+1)/4}\chi_{s_j},$$ 
where $\chi_{s_j}$, $1\le j\le\frac{q+1}{4}$, may 
not be distinct.
\end{itemize}
\end{Lemma}

\begin{Corollary}\label{char11}
Using the above notation,  
\begin{itemize}
\item[(i)] if $q\equiv 1\pmod 4$, then the character 
of $\textrm{Ind}_K^H (1_{\Cc})\cdot f_{B_s}$ is 
$\chi_s$ for each block $B_s$ of defect $0$;
\item[(ii)] if $q\equiv 3\pmod 4$, then the character 
of $\textrm{Ind}_K^H (1_{\Cc})\cdot f_{B_r}$
is $\phi_r$ for each block $B_r$ of defect $0$. 
\end{itemize}
\end{Corollary}
{\Proof} The corollary follows from Lemma~\ref{blocks} 
and Lemma~\ref{decomposition_1}.
\QED

Since $H$ preserves incidence, the following corollary 
is obvious.
\begin{Corollary}\label{u1}
Let $\p\in I$. Using the above notation, we have
$$\chi_{N_{Pa, E}(\p)}\cdot h = \chi_{N_{Pa, E}(\p^h)}, 
\chi_{N_{Se, E}(\p)}\cdot h=\chi_{N_{Se, E}(\p^h)}$$
for $h\in H$.

\end{Corollary}

In the rest of the article, we always view 
$\mathcal{G}_{\p}$ as a vector over $F$. 
Consider the maps $\phi_{\B}$, $\phi_{\D}$, and 
$\phi_{\D^{'}}$ from $F^I$ to $F^E$ defined by 
extending
$$\mathcal{G}_{\p}\mapsto \chi_{\p^\perp}, 
\mathcal{G}_{\p}\mapsto \chi_{N_{Pa, E}(\p)},
\mathcal{G}_{\p}\mapsto \chi_{N_{Se, E}(\p)}$$
linearly to $F^I$, respectively. Then it is clear 
that as $F$-linear maps, the marices of $\phi_\B$, 
$\phi_\D$, and $\phi_{\D^{'}}$ are $\B$, $\D$, and 
$\D^{'}$, respectively, and for ${\bf x}\in F^I$, 
$\phi_{\B}({\bf x})=\B {\bf x}$, 
$\phi_\D({\bf x})=\D{\bf x}$ and 
$\phi_{\D^{'}}({\bf x})=\D^{'}{\bf x}$. Moreover, 
we have the following result.
\begin{Lemma}\label{hom}
The maps $\phi_{\B}$, $\phi_{\D}$, and $\phi_{\D^{'}}$ 
are all $FH$-module homomorphisms from $F^I$ to $F^E$.

\end{Lemma}

{\Proof} Let $\mathcal{G}_{\p}$ be a basis element 
of $F^I$. Then $\phi(\mathcal{G}_{\p}\cdot h) = 
\phi(\mathcal{G}_{\p})\cdot h$ since
$$\phi_{\B}(\mathcal{G}_{\p}\cdot h)=\chi_{(\p^h)^\perp}
=\chi_{(\p^\perp)^h} = \chi_{\p^\perp}\cdot h 
=\phi_{\B}(\mathcal{G}_\p)\cdot h.$$ 
By linearity of $\phi_{\B}$, we have 
$\phi_{\B}({\bf x})\cdot h=\phi_{\B}({\bf x}\cdot h)$ 
for each ${\bf x} \in F^I$. The proof of the map 
$\phi_{\B}$ being $FH$-homomorphism is completed.

The proofs of the other two maps being homomorphisms 
are similar since $$\chi_{N_{Pa,E}(\p)}\cdot h=
\chi_{N_{Pa, E}(\p^h)}, \chi_{N_{Se, E}(\p)}\cdot h=
\chi_{N_{Se, E}(\p^h)}$$ 
for $h\in H$ and $\p\in I$ by Corollary~\ref{u1}. 
We omit the details.
\QED

For convenience, we use $\col_F(\mathbf{C})$ to denote 
the column space of the matrix $\mathbf{C}$ over $F$.

\begin{Corollary}
Using the above notation, we have 
$\Ima(\phi_\B)=\col_F(\B)$, $\Ima(\phi_\D)=\col_F(\D)$,
and $\Ima(\phi_{\D^{'}})=\col_F(\D^{'})$.
\end{Corollary}

Now we define 
$\mathcal{M}_1:=\langle\chi_{\ell}\mid \ell\in T\rangle_F$ 
and
$\mathcal{M}_2:=\langle\chi_{\ell_i}+\chi_{\ell_j}\mid \ell_i\not=\ell_j\in T\rangle_F$ 
to be the spans of the corresponding characteristic 
vectors over $F$.

\begin{Lemma}\label{u3}
The dimensions of $\mathcal{M}_1$ and $\mathcal{M}_2$ over 
$F$ are $\dim_F(\mathcal{M}_1)=q$ and $\dim_F(\mathcal{M}_2)=q-1$, 
respectively. Moreover, the all-one column vector ${\bf 1}$ 
of length $|E|$ is neither in $\mathcal{M}_1$ nor in 
$\mathcal{M}_2$.

\end{Lemma}

{\Proof} Since $\sum_{\ell \in T}\chi_{\ell} = {\bf 0}$, 
where ${\bf 0}$ is the zero column vector of $|E|$, 
it follows that $\{\chi_{\ell}\mid \ell \in T\}$
is linearly dependent over $F$, i.e. 
$\dim_F(\mathcal{M}_1)\le q$. Now let $T^{'}\subset T$ 
with $|T^{'}|=q$ and suppose that 
$\{\chi_{\ell}\mid \ell\in T^{'}\}$
is linearly dependent over $F$. Then 
$\sum_{\ell\in T^{'}}a_{\ell}\chi_{\ell}={\bf 0}$,
where $a_{\ell}\in F$ and $a_{\ell_1}\not=0$ for some 
$\ell_1\in T^{'}$. Since there are $q$ external points 
on $\ell_1$ and there are only $q-1$ tangent lines other 
than $\ell_1$ in $T^{'}$, some external point on $\ell_1$ 
must be passed only by $\ell_1$ among the tangent lines 
in $T^{'}$, which forces $a_{\ell_1}=0$, a contradiction. 
This shows that $T^{'}$ must be linearly independent 
over $F$, and so $\dim_F(\mathcal{M}_1)=q$. Moreover, 
if $T^{'}\subset T$ and $|T^{'}|=q$, then 
$\{\chi_{\ell}\mid \ell\in T^{'}\}$ must be a basis for 
$\mathcal{M}_1$.

Next if $\ell_1$ is a tangent line, then 
$\mathcal{M}_2=\langle\chi_{\ell_1}+\chi_{\ell}\mid \ell\in T\setminus\{\ell_1\}\rangle_F$ 
since $\chi_{\ell_i}+\chi_{\ell_j}=(\chi_{\ell_1}+\chi_{\ell_i})+
(\chi_{\ell_1}+\chi_{\ell_j})$. As $\sum_{\ell\in T\setminus\{\ell_1\}}(\chi_{\ell_1}+\chi_{\ell})={\bf 0}$,
$\dim_F(\mathcal{M}_2)\le q-1$. Let $T^{'}\subset T\setminus\{\ell_1\}$ 
with
$|T^{'}|=q-1$ and suppose that $\{\chi_{\ell_1}+\chi_{\ell}\mid \ell\in T^{'}\}$ 
is linearly dependent over $F$. Then 
$\sum_{\ell\in T^{'}}a_{\ell}(\chi_{\ell_1}+\chi_{\ell})
=\sum_{\ell\in T^{'}}a_{\ell}\chi_{\ell}={\bf 0}$ since 
$|T^{'}|$ is even, where $a_{\ell}\in F$ and $a_{\ell_2}\not=0$ 
for some $\ell_2\in T^{'}$. By applying the same argument 
in the first paragraph of this proof, again, we obtain that 
$a_{\ell_2}=0$ which is a contradiction. Therefore, 
$\{\chi_{\ell_1}+\chi_{\ell}\mid \ell \in T^{'}\}$ is 
linearly independent over $F$, and so $\dim_F(\mathcal{M}_2)=q-1$. 
Moreover, if $T^{'}\subset T\setminus\{\ell_1\}$ and $|T^{'}|=q-1$, 
then $\{\chi_{\ell_1}+\chi_{\ell}\mid \ell \in  T^{'}\}$ must be
a basis for $\mathcal{M}_2$.

Now we assume that ${\bf 1}\in \mathcal{M}_1$ and 
$\{\chi_{\ell}\mid \ell \in T^{'}\}$ with $T^{'}\subset T$ 
and  $|T^{'}|=q$ is a basis for $\mathcal{M}_1$. Then
$\sum_{\ell\in T^{'}}a_{\ell} \chi_{\ell}={\bf 1}$, 
where $a_{\ell}\in F$ for $\ell\in T^{'}$ and 
$a_{\ell_k}\not=0$ for some $\ell_k\in T^{'}$. Since
$|T^{'}\setminus\{\ell_k\}|=q-1$, some external point 
on $\ell_k$ must be only passed by $\ell_k$ among all 
the tangent lines in $T{'}$; this forces $a_{\ell_k}=1$. 
For each $\ell\in T'\setminus\{\ell_k\}$, we have 
$a_{\ell_k}+a_{\ell}=1$, that is, $a_{\ell}=0$ for 
each $\ell\in T'\setminus\{\ell_k\}$. Thus 
$\chi_{\ell_k}={\bf 1}$, which is impossible. 
Consequently, ${\bf 1}\notin\mathcal{M}_1$. Similarly, 
we can show that ${\bf 1}\notin \mathcal{M}_2$.
We omit the details.
\QED

\begin{Lemma}\label{u2}
If $q\equiv 1\pmod 4$, then $\col_F(\D)=\mathcal{M}_1$; 
if $q\equiv 3\pmod 4$, then $\col_F(\D)=\mathcal{M}_2$.

\end{Lemma}

{\Proof} Assume that $q\equiv 1\pmod 4$. Let 
$\chi_{N_{Pa, E}(\p)}$ be the column of $\D$ 
indexed by $\p$. Then $\chi_{N_{Pa, E}(\p)}$ 
is an $F$-linear combination of the generating 
elements of $\mathcal{M}_1$ by Corollary~\ref{sksum}. 
Now if $\chi_{\ell}$ is a generating element of 
$\mathcal{M}_1$, then it is an $F$-linear combination 
of the columns of $\D$ by Corollary~\ref{tsum1}. 
Therefore, $\col_F(\D)=\mathcal{M}_1$.

Now we assume that $q\equiv 3\pmod 4$. Let 
$\chi_{N_{Pa, E}(\p)}$ be the column of $\D$ 
indexed by $\p$. Suppose that $\ell(\p)$ is a 
tangent line through an external point on 
$\p^\perp$ and $T(\p,\ell(\p))$ is the set of 
tangent lines through the external points on 
$\ell(\p)$ that are also on the passant lines 
through $\p$. Then by Corollary~\ref{sksum} and 
the fact that $|T(\p,\ell(\p))|=\frac{q+1}{2}$ 
is even, we have
\begin{equation*}
\begin{array}{llll}
\chi_{N_{Pa,E}(\p)} & = & \displaystyle
\sum_{\ell\in T(\p, \ell(\p))}\chi_{\ell}\\
{} & =  & \displaystyle
\sum_{\ell\in T(\p,\ell(\p))}(\chi_{\ell} + 
\chi_{\ell(\p)});

\end{array}
\end{equation*}
that is, $\chi_{N_{Pa, E}(\p)}\in \mathcal{M}_2$. 
Now let $\chi_{\ell_1}+\chi_{\ell_2}$ be a 
generating element of $\mathcal{M}_2$. Then we have
\begin{equation*}
\begin{array}{lll}
\chi_{\ell_1}+\chi_{\ell_2} & = & 
\displaystyle\sum_{\q\in \mathcal{M}^{'}(\p)}
\chi_{N_{Pa, E}(\q)}
\end{array}
\end{equation*}
by Corollary~\ref{tsum2}, where $\p=\ell_1\cap \ell_2$. 
Hence, $\col_F(\D)=\mathcal{M}_2$.
\QED

\begin{Corollary}\label{dim}
If $q\equiv 1\pmod 4$, $\rank_2(\D)=q$; if 
$q\equiv 3\pmod 4$, $\rank_2(\D)=q-1$.

\end{Corollary}
{\Proof} It follows from Lemmas~\ref{u3} and 
\ref{u2}. \QED

Further, we have the following decomposition 
of $\col_F(\D^{'})$.

\begin{Lemma}\label{deofD}
If $q\equiv 3\pmod 4$, then $\col_F(\D^{'})=
\langle{\bf 1}\rangle\oplus \col_F(\D)$
as $FH$-modules, where $\langle{\bf 1}\rangle$ 
is the trivial $FH$-module generated by the 
all-one column vector ${\bf 1}$.

\end{Lemma}

{\Proof } Since each row of $\D^{'}$ has 
$\frac{(q-1)^2}{4}$ $1$s, then
$$\displaystyle\sum_{\p\in I}\chi_{N_{Se, E}(\p)}
={\bf 1}.$$ For $h\in H$,
\begin{equation*}
{\bf 1}\cdot h=(\displaystyle\sum_{\p\in I}\chi_{N_{Se, E}(\p)})\cdot h
=\displaystyle\sum_{\p\in I}\chi_{N_{Se, E}(\p^h)}
=\displaystyle\sum_{\p\in I}\chi_{N_{Se, E}(\p)}
={\bf 1}\in \col_F(\D).
\end{equation*}
Consequently, $\langle{\bf 1}\rangle$ is indeed a 
trivial submodule of $\col_F(\D^{'})$.

It is clear that 
$\col_F(\D^{'})=\langle{\bf 1}\rangle + \col_F(\D)$ 
since $\chi_{N_{Se, E}(\p)}\in \col_F(\D^{'})$ 
if and only if 
$\chi_{N_{Se, E}(\p)}={\bf 1}+\chi_{N_{Pa, E}(\p)}\in \langle{\bf 1}\rangle
+\col_F(\D)$. Further, 
$\langle{\bf 1}\rangle\cap \col_F(\D)={\bf 0}$ since 
$\col_F(\D)=\mathcal{M}_2$ and ${\bf 1}\notin \mathcal{M}_2$ 
by Lemmas ~\ref{u3} and ~\ref{u2}. Therefore, 
$\col_F(\D^{'})=\langle{\bf 1}\rangle\oplus \col_F(\D)$.
\QED

\section{Statement and Proof of Main Theorem}
The main theorem is given as follows.

\begin{Theorem}\label{main}
Let $\Ima(\phi_\B)$ and $\Ima(\phi_\D)$ be defined 
as above. As $FH$-modules,
\begin{itemize}
\item[(i)] if $q\equiv 1\pmod 4$, then

$$\Ima(\phi_\B)=\Ima(\phi_\D)\oplus (\displaystyle\bigoplus_{s=1}^{(q-1)/4}M_s),$$
where $M_s$ for $1\le s\le \frac{q-1}{4}$ are 
pairwise non-isomorphic simple $FH$-modules of 
dimension $q-1$;

\item[(ii)] if $q\equiv 3\pmod 4$, then

$$\Ima(\phi_\B)=\langle{\bf 1}\rangle\oplus \Ima(\phi_\D)\oplus (\displaystyle\bigoplus_{r=1}^{(q-3)/4}M_r),$$
where $M_r$ for $1\le s\le \frac{q-3}{4}$ are pairwise 
non-isomorphic simple $FH$-modules of dimension $q+1$ 
and $\langle{\bf 1}\rangle$ is the trivial $FH$-module 
generated by the all-one column vector of length $|E|$.

\end{itemize}

\end{Theorem}
To prove the main theorem, we need the following lemma.

\begin{Lemma}\label{y4}
Let  $q-1=2^n m$ or $q+1=2^n m$
with $2\nmid m$ depending on whether $q\equiv 1\pmod 4$ 
or $q\equiv 3\pmod 4$. Using the above notation, 
\begin{itemize}
\item[(i)] if $q\equiv 1\pmod 4$, then 
$\Ima(\phi_\B)\cdot e_{B_0}=\Ima(\phi_\D)$,
$\Ima(\phi_\B)\cdot e_{B_s} \not={\bf 0}$ 
for $1\le s\le \frac{q-1}{4}$, and 
$\Ima(\phi_\B)\cdot e_{B_t^{'}}={\bf 0}$ for 
$m\ge 3$ and $1\le t \le \frac{m-1}{2}$;

\item[(ii)] if $q\equiv 3\pmod 4$, then 
$\Ima(\phi_\B)\cdot e_{B_0}=\Ima(\phi_{\D^{'}})$, 
$\Ima(\phi_\B)\cdot e_{B_r} \not={\bf 0}$ for 
$1 \le r \le \frac{q-3}{4}$, and 
$\Ima(\phi_\B)\cdot e_{B_t^{'}}={\bf 0}$ for 
$m\ge 3$ and $1\le t\le \frac{m-1}{2}$.
\end{itemize}

\end{Lemma}

{\Proof} It is clear that $\Ima(\phi_\B)$ is 
generated by $\{\chi_{\p^\perp}\mid \p\in I\}$
over $F$. Let $B\in Bl(H)$. Since
\begin{equation*}
\begin{array}{llll}
{\chi_{\p^\perp}}\cdot e_B &= &\displaystyle
\sum_{C\in cl(H)}e_{B}(\widehat{C})\displaystyle
\sum_{h\in C}{\chi}_{\p^\perp}\cdot h\\
{} & = & \displaystyle\sum_{C\in cl(H)}e_{B}(\widehat{C})
\displaystyle\sum_{h\in C}{\chi}_{(\p^\perp)^h},\\
{} & = & \displaystyle\sum_{C\in cl(H)}e_{B}(\widehat{C})
\displaystyle\sum_{h\in C}\sum_{\q\in(\p^{\perp})^h\cap E}
{\chi}_{\q},

\end{array}
\end{equation*}
we have
\begin{equation*}
{\chi_{\p^\perp}}\cdot e_B =\sum_{\q\in I}\mathcal{S}
(B,\p,\q){\chi}_{\q},
\end{equation*}
where
\begin{equation*}
\begin{array}{lll}

S(B,\p,\q):=\displaystyle\sum_{C\in cl(H)}
|\mathcal{H}_{\p,\q}\cap C|e_B(\widehat{C}).
\end{array}
\end{equation*}

Assume first that $q\equiv 1\pmod 4$. If 
$\ell_{\p,\q}\in Se_{\p}$, then $S(B,\p,\q)=0$ 
for each $B\in Bl(H)$ since $|\mathcal{H}_{\p,\q}\cap C|=0$
in $F$ for each $C\not=[0]$ by Lemma~\ref{m1}(i) 
and 
$e_{B_0}(\widehat{[0]})=e_{B_s}(\widehat{[0]})=e_{B_t^{'}}(\widehat{[0]})=0$ 
by 1(c), 2(c), 3(c) of Lemma~\ref{expression}.

If $\ell_{\p,\q}\in Pa_\p$ and $\q\in \p^\perp$, then by 
Lemma~\ref{m1}(ii) and 1(a), 1(c), 2(a), 2(c), 3(a), 3(c)
of Lemma~\ref{expression},
\begin{equation*}
\begin{array}{llllllll}
S(B_0,\p,\q)& = & |\mathcal{H}_{\p,\q}\cap [0]| e_{B_0}(\widehat{[0]})+|\mathcal{H}_{\p,\q}\cap D|
e_{B_0}(\widehat{D}) & = & 0+1 & = &1,\\
S(B_s,\p,\q)& = & |\mathcal{H}_{\p,\q}\cap [0]| e_{B_s}(\widehat{[0]})+|\mathcal{H}_{\p,\q}\cap D|
e_{B_s}(\widehat{D}) & = & 0+0 & = &0,\\
S(B_t^{'},\p,\q)& = & |\mathcal{H}_{\p,\q}\cap [0]| e_{B_t^{'}}(\widehat{[0]})+|\mathcal{H}_{\p,\q}\cap D|
e_{B_t^{'}}(\widehat{D}) & = & 0+0 & =& 0.\\
\end{array}
\end{equation*}

If $\q$ is on a passant line $\ell$ through $\p$
and $\q\notin\p^\perp$, 
then by Lemma~\ref{m1}(iii) and
1(c), 1(d), 2(c), 2(d), 3(c), 3(d) of Lemma~\ref{expression}, 
\begin{equation*}
\begin{array}{llllllll}
S(B_0,\p,\q)& = & |\mathcal{H}_{\p,\q}\cap [0]| e_{B_0}(\widehat{[0]})+|\mathcal{H}_{\p,\q}\cap [\pi_k]|
e_{B_0}(\widehat{[\pi_k]}) & = & 1, & {}\\
S(B_s,\p,\q)& = & |\mathcal{H}_{\p,\q}\cap [0]| e_{B_s}(\widehat{[0]})+|\mathcal{H}_{\p,\q}\cap [\pi_k]|
e_{B_s}(\widehat{[\pi_k]}) & = & e_{B_s}(\widehat{[\pi_k]}), & {}\\
S(B_t^{'},\p,\q)& = & |\mathcal{H}_{\p,\q}\cap [0]| e_{B_t^{'}}(\widehat{[0]})+|\mathcal{H}_{\p,\q}\cap [\pi_k]|
e_{B_t^{'}}(\widehat{[\pi_k]}) & = & 0. & {}\\
\end{array}
\end{equation*}
By Lemma~\ref{m1}(iii) and the fact
that there are $\frac{q-1}{4}$ classes
of the form $[\pi_k]$ and there are 
$\frac{q-1}{2}$ points on $\ell$
that are not on $\p^\perp$, we have that
for each $[\pi_k]$ there exist two external
points $\q_1$ and $\q_2$ on $\ell$ such 
that $|\mathcal{H}_{\p,\q_j}\cap [\pi_k]|$ 
($j=1$ or $2$) is odd and for each 
$\q\in\ell$ and $\q\notin\p^\perp$ 
there is a class $[\pi_k]$ such that
$|\mathcal{H}_{\p,\q}\cap [\pi_k]|$
is odd. Combining the above analysis with
Lemma~\ref{EX}, we obtain that for 
each $B_s$, there is a $\q$ and a 
class $[\pi_k]$ such that 
$S(B_s,\p,\q)=e_{B_s}(\widehat{[\pi_k]})\not=0$.

Therefore, we have shown that 
$\Ima(\phi_\B)\cdot e_{B_0}=\Ima(\phi_\D)$ by 
definition, $\Ima(\phi_\B)\cdot e_{B_s}\not={\bf 0}$ 
for each $s$, and 
$\Ima(\phi_{\B})\cdot e_{B_t^{'}}={\bf 0}$.
The proof of (i) is completed.

Now assume that $q\equiv 3\pmod 4$. 
If $\ell_{\p,\q}\in Pa_{\p}$, then 
$S(B,\p,\q)=0$ for each $B\in Bl(H)$ 
since $|\mathcal{H}_{\p,\q}\cap C|=0$ 
by 4(d), 5(c), 6(d) of Lemma~\ref{expression}.

Let $\ell_{\p,\q}\in Se_\p$ and $\q\in \p^\perp$, 
then by Lemma~\ref{m2}(ii) and 4(a), 4(d), 5(a), 
5(c), 5(d) of Lemma~\ref{expression}
\begin{equation*}
\begin{array}{llllllllll}
S(B_0,\p,\q)& = & |\mathcal{H}_{\p,\q}\cap [0]| e_{B_0}(\widehat{[0]})+|\mathcal{H}_{\p,\q}\cap D|
e_{B_0}(\widehat{D}) & = & 0+1 & =&1,\\
S(B_s,\p,\q)& = & |\mathcal{H}_{\p,\q}\cap [0]| e_{B_s}(\widehat{[0]})+|\mathcal{H}_{\p,\q}\cap D|
e_{B_s}(\widehat{D}) & = & 0+0 & =& 0,\\
S(B_t^{'},\p,\q)& = & |\mathcal{H}_{\p,\q}\cap [0]| e_{B_t^{'}}(\widehat{[0]})+|\mathcal{H}_{\p,\q}\cap D|
e_{B_t^{'}}(\widehat{D}) & = & 0+0 &=& 0.\\
\end{array}
\end{equation*}

If $\ell_{\p,\q}\in Se_{\p}$ and $\q\notin\p^\perp$, 
then by Lemma~\ref{m2}(i), 4(c), 4(d), 5(c), 5(d), 
6(c), 6(d) of Lemma~\ref{expression},
\begin{equation*}
\begin{array}{llllllll}
S(B_0,\p,\q)& = & |\mathcal{H}_{\p,\q}\cap [0]| e_{B_0}(\widehat{[0]})+|\mathcal{H}_{\p,\q}\cap [\theta_i]|
e_{B_0}(\widehat{[\theta_i]}) & = & 1, & {}\\
S(B_s,\p,\q)& = & |\mathcal{H}_{\p,\q}\cap [0]| e_{B_s}(\widehat{[0]})+|\mathcal{H}_{\p,\q}\cap [\theta_i]|
e_{B_s}(\widehat{[\pi_k]}) & = & e_{B_s}(\widehat{[\theta_i]}), & {}\\
S(B_t^{'},\p,\q)& = & |\mathcal{H}_{\p,\q}\cap [0]| e_{B_t^{'}}(\widehat{[0]})+|\mathcal{H}_{\p,\q}\cap [\theta_i]|
e_{B_t^{'}}(\widehat{[\theta_i]}) & = & 0. & {}\\
\end{array}
\end{equation*}
From Lemma~\ref{m2}(i) and Lemma~\ref{EX}, we have 
that for each $B_s$, there is a $\q$ and a class 
$[\theta_i]$ such that 
$S(B_s, \p,\q)=e_{B_s}(\widehat{[\theta_i]})\not=0$.

Therefore, we have shown that 
$\Ima(\phi_\B)\cdot e_{B_0}=\Ima(\phi_{\D^{'}})$ 
by definition, $\Ima(\phi_{\B})\cdot e_{B_s}\not={\bf 0}$ 
for each $s$, and $\Ima(\phi_\B)\cdot e_{B_t^{'}}={\bf 0}$.
The proof of (ii) is completed.
\QED

\noindent{{\bf Proof of Theorem~\ref{main}:}}  
Let $B$ be a $2$-block of defect $0$ of $H$. 
Then by Lemma~\ref{reduction}, we have
$$F^I\cdot e_B=\overline{\mathbf{S}^I\cdot f_B}.$$
Therefore, by Corollary~\ref{char11}, $F^I\cdot e_b=N$, 
where $N$ is the simple $FH$-module of dimension 
$q-1$ or $q+1$ lying in $B$ accordingly as 
$q\equiv 1\pmod 4$ or $q\equiv 3\pmod 4$. It is clear 
that $\phi_\B(F^I)=\Ima(\phi_\B)$.

Assume that $q\equiv 1\pmod 4$ and $q-1=m2^n$ with 
$2\nmid m$. Since
$$1=e_{B_0}+\displaystyle\sum_{s=1}^{(q-1)/4}e_{B_s}
+\displaystyle\sum_{t=1}^{(m-1)/2}e_{B_t^{'}},$$
we have
\begin{equation}\label{last1}
\begin{array}{lllllll}
\Ima(\phi_\B) & = & \Ima(\phi_\B)\cdot e_{B_0}\oplus (\displaystyle\bigoplus_{s=1}^{(q-1)/4}\Ima(\phi_\B)\cdot e_{B_s})
\oplus (\displaystyle\bigoplus_{t=1}^{(m-1)/2}\Ima(\phi_\B)\cdot e_{B_t^{'}})\\
{} & = & \Ima(\phi_\D)\oplus (\displaystyle\bigoplus_{s=1}^{(q-1)/4}\phi_{\B}(F^I)\cdot e_{B_s})\\
{} & = & \Ima(\phi_\D)\oplus (\displaystyle\bigoplus_{s=1}^{(q-1)/4}\phi_{\B}(F^I\cdot e_{B_s}))\\
{} & = & \Ima(\phi_\D)\oplus (\displaystyle\bigoplus_{s=1}^{(q-1)/4}\phi_{\B}(N_s))\\
{} & = & \Ima(\phi_\D)\oplus (\displaystyle\bigoplus_{s=1}^{(q-1)/4}M_s),

\end{array}
\end{equation}
where $N_s$ is the simple module of dimension $q-1$ 
lying in $B_s$ for each $s$ by the discussion in the 
first paragraph and $M_s:=\phi_\B(N_s)$ for each $s$. 
In (\ref{last1}), the terms $e_{B_t^{'}}$ for 
$1\le t\le \frac{m-1}{2}$ and $\Ima(\phi_{\B})\cdot e_{B_t^{'}}$
for $1\le t\le \frac{m-1}{2}$ appear only when $m\ge 3$; 
the second equality holds since 
$\Ima(\phi_\B)\cdot e_{B_t^{'}}={\bf 0}$ for each $t$ 
and $\Ima(\phi_\B)\cdot e_{B_0}= \Ima(\phi_\D)$
by Lemma~\ref{y4}(i); and the third equality holds 
since $\phi_\B$ is an $FH$-homomorphism by Lemma~\ref{hom} 
and $e_{B_s}\in FH$. Consider the map
$$\lambda_S: N_S\rightarrow \phi_\B(N_s)$$
defined by $\lambda_s(n)=\phi_\B(n)$ for $n\in N_s$, 
where $1\le s\le \frac{q-1}{4}$. It is clear that 
$\lambda_s$ is the same as the resctriction of 
$\phi_\B$ to $N_s$. Consequently, $\lambda_s$ is a 
surjective $FH$-homomorphism. Moreover, $\Ker(\lambda_s)$ 
is either ${\bf 0}$ or $N_s$ since, otherwise, 
$\Ker(\lambda_s)$ would be a non-trivial submodule of $N_s$
which is impossible. If $\Ker(\lambda_s)=N_s$, then 
$\phi_\B(N_s)=\phi_{\B}(F^I)\cdot e_{B_s}={\bf 0}$,
which is not the case by Lemma~\ref{y4}(i). Thus, 
we must have $\Ker(\lambda_s)={\bf 0}$; that is,
$\lambda_s$ is an $FH$-isomorphism. So we have shown 
that $M_s:=\Ima(N_s)\cong N_s$ and thus $M_s$ for 
$1\le s\le \frac{q-1}{4}$ are pairwise non-isomorphic 
simple modules of dimension $q-1$. The proof of (i) 
is finished.

Now assume that $q\equiv 3\pmod 4$. Applying the same 
argument as above, we have
$$\Ima(\phi_\B)=\Ima(\phi_{\D^{'}})\oplus (\displaystyle\bigoplus_{r=1}^{(q-3)/4}M_r),$$
where $M_r$ for $1\le r\le \frac{q-3}{4}$ are pairwise 
non-isomorphic simple $FH$-modules of dimension $q+1$.
Since $\Ima(\phi_{\D^{'}})=\langle{\bf 1}\rangle\oplus \Ima(\phi_\D)$ 
by Lemma~\ref{deofD}, it follows that
$$\Ima(\phi_\B)=\langle {\bf 1}\rangle\oplus \Ima(\phi_\D)\oplus(\displaystyle\bigoplus_{r=1}^{(q-3)/4}M_r).$$
\QED

Now Conjecture~\ref{conj} follows as a corollary.

\begin{Corollary}
Let $\mathcal{L}$ and $\mathcal{L}_0$ be the 
$\Ff_2$-null spaces of $\B$ and $\B_0$, respectively. 
Then
\begin{equation*}
\dim_{\Ff_2}(\mathcal{L})=
\begin{cases}
\frac{q^2-1}{4}-q, & \text{if}\;q\equiv 1\pmod 4,\\
\frac{q^2-1}{4}-q +1,& \text{if}\; q\equiv 3\pmod 4\\
\end{cases}
\end{equation*}
and

\begin{equation*}
\dim_{\Ff_2}(\mathcal{L}_0)=
\begin{cases}
\frac{q^2-1}{4}, & \text{if}\;q\equiv 1\pmod 4,\\
\frac{q^2-1}{4}+1,& \text{if}\; q\equiv 3\pmod 4.\\
\end{cases}
\end{equation*}
\end{Corollary}
{\Proof} From Theorem~\ref{main} and Corollary~\ref{dim}, 
it follows that the $2$-rank of
$\B$ is 
$$\rank_2(\B)=q+\frac{(q-1)^2}{4}$$
or
$$\rank_2(\B)=1+(q-1)+\frac{(q-1)^2}{4}$$
accordingly as $q\equiv 1\pmod 4$ or $q\equiv 3\pmod 4$. 
Therefore, the dimension of the $\Ff_2$-null space of 
$\B$ is
$$\dim_{\Ff_2}(\mathcal{L})=\frac{q(q-1)}{2}-(q+\frac{(q-1)^2}{4})=\frac{q^2-1}{4}-q$$
or
$$\dim_{\Ff_2}(\mathcal{L})=\frac{q(q-1)}{2}-(1+(q-1)+\frac{(q-1)(q-3)}{4})=\frac{q^2-1}{4}-q+1$$
accordingly as $q\equiv 1\pmod 4$ or $q\equiv 3\pmod 4$.

Since $\rank_2(\B)=\rank_2(\B_0)$, the dimension of 
$\mathcal{L}_0$ can be calculated in the same way. 
We omit the details.
\QED

\end{document}